        %%%%%%%%%%%%%%%%%%%%%%%%%%%%%%%%%%%%%%%
        %%%                                 %%%
        %%     QUANTUM FUNCTION ALGEBRAS     %%
        %%                 AS                %%
        %%    QUANTUM ENVELOPING ALGEBRAS    %%
        %%                                   %%
        %%          by Fabio Gavarini        %%
        %%%                                 %%%
        %%%%%%%%%%%%%%%%%%%%%%%%%%%%%%%%%%%%%%%

%&amstex
\input amstex
\documentstyle{amsppt}

\magnification=\magstep1
\hsize=6.5truein
\vsize=9truein

\font \smallrm=cmr10 at 10truept
\font \smallbf=cmbx10 at 10truept
\font \smallit=cmti10 at 10truept
 at 10truept
 at 10truept
\font \smallsl=cmsl10 at 10truept

\baselineskip=.15truein

\font\cs=eufm10

\def \llongrightarrow
{\,\relbar\joinrel\relbar\joinrel\relbar\joinrel\rightarrow\,}
\def \llonghookrightarrow
{\,\lhook\joinrel\relbar\joinrel\relbar\joinrel\relbar\joinrel\rightarrow\,}
\def \llongtwoheadrightarrow
{\,\relbar\joinrel\relbar\joinrel\relbar\joinrel\twoheadrightarrow\,}
\def \N {{\Bbb N}}
\def \kq {k(q)}
\def \qm {q^{-1}}
\def \kqqm {k \! \left[ q,\qm \right]}
\def \gerg {\hbox{{\cs g}}}
\def \gerh {\hbox{{\cs h}}}
\def \gerb {\hbox{{\cs b}}}
\def \gern {\hbox{{\cs n}}}
\def \gerU {\hbox{{\cs U}}}
\def \gerF {\hbox{{\cs F}}}
\def \calU {\hbox{$ \Cal U $}}
\def \calF {\hbox{$ \Cal F $}}
\def \slgot {\hbox{{\cs sl}}}
\def \gl {\hbox{{\cs gl}}}
\def \fbar {\overline F}
\def \ebar {\overline E}
\def \ug {U(\gerg)}
\def \uh {U(\gerh)}

\def \uqsln {U_q \big( \slgot(n+1) \big)}
\def \uqQsln {U_q^{\scriptscriptstyle Q} \big( \slgot(n+1) \big)}
\def \uqPsln {U_q^{\scriptscriptstyle P} \big( \slgot(n+1) \big)}
\def \uqgln {U_q \big( \gl(n+1) \big)}

\def \fqPg {F_q^{\scriptscriptstyle P}[G]}
\def \fqPsln {F_q^{\scriptscriptstyle P}[SL(n+1)]}
\def \fqQsln {F_q^{\scriptscriptstyle Q}[SL(n+1)]}
\def \fqQtildesln {\widetilde{F}_q^{\scriptscriptstyle Q}[SL(n+1)]}
\def \fqPtildesln {\widetilde{F}_q^{\scriptscriptstyle P}[SL(n+1)]}
\def \fqtildesln {\widetilde{F}_q[SL(n+1)]}
\def \funoQtildesln {\widetilde{F}_1^{\scriptscriptstyle Q}[SL(n+1)]}
\def \funoPtildesln {\widetilde{F}_1^{\scriptscriptstyle P}[SL(n+1)]}
\def \funotildesln {\widetilde{F}_1[SL(n+1)]}
\def \calfQsln {\calF^{\scriptscriptstyle Q}[SL(n+1)]}
\def \calfPsln {\calF^{\scriptscriptstyle P}[SL(n+1)]}
\def \fqmn {F_q[M(n+1)]}

%\NoBlackBoxes

\document

\topmatter

{\ }

\vskip-33pt

\hfill   {{\smallsl Communications in Algebra\/}
{\smallbf 26},  {\smallrm no.~6 (1998), 1795--1818}}
\hskip19pt   {\ }

\vskip41pt

\title
   Quantum function algebras as quantum enveloping algebras
\endtitle

\author
   Fabio Gavarini
\endauthor

\affil
   Institut de Recherche Math\'ematique Avanc\'ee, Universit\'e{} "Louis
Pasteur" -- C.N.R.S.  \\
   7, rue Ren\'e{} Descartes, 67084 Strasbourg Cedex, France  \\
\endaffil

\address\hskip-\parindent
        Institut de Recherche Math\'ematique Avanc\'ee  \newline
        Universit\'e{} "Louis Pasteur" -- C.N.R.S.  \newline
        7, rue Ren\'e{} Descartes  \newline
        67084 STRASBOURG Cedex --- FRANCE  \newline
        e-mail: \  gavarini\@math.u-strasbg.fr  \newline
        \phantom{e-mail:} \, \  gavarini\@mat.uniroma1.it  \newline
        \phantom{e-mail:} \, \  gavarini\@mat.uniroma3.it
\endaddress

\abstract
   Inspired by a result in [Ga], we locate three integer forms of
$ F_q[SL(n+1)] $  over  $ \kqqm $,  with a presentation by generators
and relations, which for  $ \, q=1 \, $  specialize to  $ U(\gerh) $,
where  $ \gerh $  is the Lie bialgebra of the Poisson Lie group dual
to  $ SL(n+1) $.  In sight of this we prove two PBW-like theorems for
$ F_q[SL(n+1)] $,  both related to the classical PBW theorem for
$ U(\gerh) $.
\endabstract

\endtopmatter

\footnote""{ 1991 {\it Mathematics Subject Classification:}
Primary 17B37, 81R50 }
\footnote""{ Partially supported by a post-doc fellowship of the
{\it Consiglio Nazionale delle Ricerche} \, (Italy) }

\vskip0,9truecm

\centerline { \bf  Introduction }

\vskip12pt

\hfill  \hbox{\vbox{
     \hbox{\it \hskip10pt  "Nel mezzo a una  $q$--algebra  di funzioni }
     \hbox{\it \hskip26pt    io ci ritrovo una tal forma intera }
     \hbox{\it \hskip26pt        che l'algebra di Lie dual mi doni" }
        \vskip4pt
     \hbox{\sl \hskip47pt     N.~Barbecue, "Scholia" } }
\hskip1truecm }

\vskip8pt

  Let  $ G $  be a connected, simply connected, semisimple
algebraic group over an algebraically closed field  $ k $  of characteristic
zero, and consider on it the
Sklyanin-Drinfel'd structure of Poisson group (cf.~for instance [DP] \S 11 or
[Ga] \S 1, or even [Dr]); then  $ \, \gerg := Lie(G) \, $  is a Lie bialgebra,
$ F[G] $  is a Poisson Hopf algebra, and  $ U(\gerg) $  is a Poisson Hopf
coalgebra.  Let  $ H $  be the corresponding dual Poisson
(algebraic) group of  $ G $,  whose tangent Lie bialgebra
$ \, \gerh := Lie(H) \, $  is the (linear) dual of  $ \gerg \, $:  then again
$ F[H] $  is a Poisson Hopf algebra, and
$ U(\gerh) $  is a Poisson Hopf coalgebra.
                                                    \par
  The quantum group  $ U_q^{\scriptscriptstyle Q}(\gerg) $  of Drinfel'd and
Jimbo provides a quantization of
$ U(\gerg) $:  namely,  $ U_q^{\scriptscriptstyle Q}
(\gerg) $  is a Hopf algebra over  $ \kq $  which has a
$ \kqqm $--form  $ \gerU^{\scriptscriptstyle Q}(\gerg) $  which for  $ \, q
\rightarrow 1 \, $  specializes to  $ \ug $  as a Poisson Hopf coalgebra.
Dually, by means of a Peter-Weyl type axiomatic trick one constructs a Hopf
algebra  $ F_q^{\scriptscriptstyle P}[G] $  of matrix
coefficients of  $ U_q^{\scriptscriptstyle Q}(\gerg) $  with a  $ \kqqm $--form
$ \gerF^{\scriptscriptstyle P}[G] $  which specializes to  $ F[G] $,  as a
Poisson Hopf algebra,
for  $ \, q \rightarrow 1 \, $.  So far the quantization only dealt with the
Poisson group
$ G \, $;  the dual group  $ H $  is involved defining a different  $ \kqqm
$--form  $ \calU^{\scriptscriptstyle P}(\gerg) $  (of a quantum group  $
U_q^{\scriptscriptstyle P}(\gerg) $)  which specializes to  $ F[H] $  (as a
Poisson Hopf algebra) for  $ \, q \rightarrow 1 \, $  (cf.~[DP] or [DKP]).  In a
dual fashion, it is proved in [Ga]   --- in a wider context ---   that the dual
(in the Hopf sense) quantum function algebra
$ F_q^{\scriptscriptstyle Q}[G] $  has a  $ \kqqm $--integer  form  $
\calF^{\scriptscriptstyle Q} [G] $  which for  $ \, q \rightarrow 1 \, $
specializes to  $ U(\gerh) $,  as a Poisson Hopf coalgebra.  Therefore quantum
function algebras can also be thought of as quantum enveloping algebras, whence
the title of the paper.
                                                 \par
  In this paper we stick to the case of the group  $ \, G = SL(n+1) $.
                                                 \par
  Our first goal is to relate the latter result above with the well-known
presentation of  $ \fqPsln $  by generators and relations (cf.~[FRT]): namely,
inspired by the definition of  $ \calF^{\scriptscriptstyle Q}[G] $  and
$ \calF^{\scriptscriptstyle P}[G] $,  we define two
$ \kqqm $--integer  forms  $ \fqQtildesln $  and
$ \fqPtildesln $  (along with a third one,  $ \fqtildesln $)  of  $ \fqPsln $;
these inherit a presentation by generators and relations, which enables us to
prove that they specialize to  $ U(\gerh) $  (as a Poisson Hopf coalgebra) for
$ \, q \rightarrow 1 $.  As a second step, since for  $ U(\gerh) $  one has the
Poincar\'e-Birkhof-Witt (PBW in short) theorem which provides "monomial" basis,
because of the previous result we are led to look for PBW-like theorems for
$ \fqPsln $:  we provide two of them, both closely related with the classical
PBW theorem for  $ U(\gerh) $.
                                                 \par
  The paper is organized as follows.  Sections 1, 2 are introductory.  Sections
3, 4 are devoted to integer forms of
$ \fqPsln $  and their specialization.  Section 5 is an excursus, where we
explain the relation among the constructions and results in this paper and those
in [Ga]: this is one of the main motivation of this work; on the other hand,
this section can be skipped without affecting the comprehension of the rest of
the paper, which is completely self-contained.  In section 6 we briefly outline
the extension of the previous results to the quantum function algebra
$ F_q[GL(n+1)] $.  Finally, section 7 deals with PBW theorems.

\vskip1,3truecm

\centerline { \bf  \S \; 1 \,  The universal enveloping algebra
$ \uh $ }

\vskip10pt

  A presentation of  $ U(\gerh) $  by generators and relations (in the general
case) can be found in [Ga], \S 1.  When  $ \, G = SL(n+1) \, $,  it reads as
follows.
                                                 \par
  $ U(\gerh) $  is the associative  $ k $--algebra  with 1 generated by
$ \text{f}_1, \dots, \text{f}_n $,  $ \text{h}_1, \dots, \text{h}_n $,
$ \text{e}_1, \dots, \text{e}_n $  (which are to be thought of as "Chevalley
generators")  with relations
  $$  \eqalignno{
   {}  &  \text{h}_i \text{h}_j - \text{h}_j \text{h}_i = 0  \hskip50pt  &
\forall\; i, j  \qquad  \cr
   \text{h}_i \text{f}_j - \text{f}_j \text{h}_i  &  = \big( 2 \, \delta_{i,j} -
\delta_{i-1,j} - \delta_{i+1,j} \big) \, \text{f}_j  \hskip50pt  &   \forall\;
i, j  \qquad  \cr
   \text{h}_i \text{e}_j - \text{e}_j \text{h}_i  &  = \big( 2 \, \delta_{i,j} -
\delta_{i-1,j} - \delta_{i+1,j} \big) \, \text{e}_j  \hskip50pt  &   \forall\;
i, j  \qquad  \cr
   {}  &  \text{f}_i \text{f}_j - \text{f}_j \text{f}_i = 0  \hskip50pt  &
\forall\; i, j :  \vert i-j \vert > 1  \qquad  \cr
   {}  &  \text{e}_i \text{e}_j - \text{e}_j \text{e}_i = 0  \hskip50pt  &
\forall\; i, j :  \vert i-j \vert > 1  \qquad  \cr
     \text{f}_i^{\,2} \text{f}_j  &  - 2 \, \text{f}_i \text{f}_j \text{f}_i +
\text{f}_j \text{f}_i^{\,2} = 0  \hskip50pt  &   \forall\; \vert i-j \vert = 1
\qquad  \cr
     \text{e}_i^{\,2} \text{e}_j  &  - 2 \, \text{e}_i \text{e}_j \text{e}_i +
\text{e}_j \text{e}_i^{\,2} = 0  \hskip50pt  &   \forall\; \vert i-j \vert = 1
\qquad  \cr
     {}  &  \text{f}_i \text{e}_j - \text{e}_j \text{f}_i = 0  \hskip50pt  &
\forall\; i, j  \qquad  \cr }  $$
   \indent   Furthermore  $ U(\gerh) $  has a Poisson Hopf coalgebra structure
given by
  $$  \eqalign{
   \Delta(\text{f}_i) = \text{f}_i \otimes 1 + 1 \otimes \text{f}_i \, ,  \qquad
S(\text{f}_i)  &  = -\text{f}_i \, ,  \qquad  \epsilon(\text{f}_i) = 0  \cr
   \Delta(\text{h}_i) = \text{h}_i \otimes 1 + 1 \otimes \text{h}_i \, ,  \qquad
S(\text{h}_i)  &  = -\text{h}_i \, ,  \qquad  \epsilon(\text{h}_i) = 0  \cr
   \Delta(\text{e}_i) = \text{e}_i \otimes 1 + 1 \otimes \text{e}_i \, ,  \qquad
S(\text{e}_i)  &  = -\text{e}_i \, ,  \qquad  \epsilon(\text{e}_i) = 0
\cr }  $$
for all  $ \, i = 1, \dots, n \, $,  and by
  $$  \displaylines{
   \delta(\text{f}_i) = \text{h}_i \wedge \text{f}_i + 2 \cdot \Bigg(
\sum_{j=1}^{i-1} \text{f}_{i+1,j} \wedge \text{e}_{j,i} + \sum_{j=i+2}^{n+1}
\text{e}_{i+1,j} \wedge \text{f}_{j,i} \Bigg)  \cr
   \delta(\text{h}_i) = 4 \cdot \Bigg( \sum_{j=1}^{i-1} \text{f}_{i,j} \wedge
\text{e}_{j,i} + \sum_{j=i+1}^{n+1} \text{e}_{i,j} \wedge \text{f}_{j,i} -
\sum_{j=1}^i \text{f}_{i+1,j} \wedge \text{e}_{j,i+1} - \sum_{j=i+2}^{n+1}
\text{e}_{i+1,j} \wedge \text{f}_{j,i+1} \Bigg)  \cr
   \delta(\text{e}_i) = \text{e}_i \wedge \text{h}_i + 2 \cdot \Bigg(
\sum_{j=1}^{i-1} \text{e}_{j,i+1} \wedge \text{f}_{i,j} + \sum_{j=i+2}^{n+1}
\text{f}_{j,i+1} \wedge \text{e}_{i,j} \Bigg)  \cr }  $$
for all  $ \, i = 1, \dots, n \, $,  where  $ \, x \wedge y := x \otimes y - y
\otimes x \, $  and the symbols  $ \, \text{f}_{h,k} \, $,  $ \, \text{e}_{h,k}
\, $,  have the following meaning:
  $$  \hbox{ $ \eqalign{
   \text{e}_{i,i+1} := \text{e}_i \; ,  &  \qquad  \text{e}_{i,j} := - \left[
\text{e}_{i,j-1}, \text{e}_{j-1,j} \right] = \left[ \text{e}_{j-1,j},
\text{e}_{i,j-1} \right]  \qquad \forall \; i < j-1  \cr
   \text{f}_{j+1,j} := \text{f}_j \; ,  &  \qquad  \text{f}_{j,i} := \left[
\text{f}_{j-1,i}, \text{f}_{j,j-1} \right] = - \left[ \text{f}_{j,j-1},
\text{f}_{j-1,i} \right]  \qquad \forall \; j > i+1  \cr } $ }   \eqno (1.1)  $$
the symbol  $ \, [\ ,\ ] \, $  denoting the usual commutator.  In fact, if
$ M_{i,j} $  ($ i, j \in \{\,1,2,\dots,n+1\,\} $)  denotes the square matrix of
size  $ n+1 $  with a 1 as  $ (i,j) $--th  entry  and all other entries equal to
 $ 0 $,  the recipe  $ \, \text{e}_h \mapsto M_{h,h+1} \, \forall\, h=1,\dots,n
\, $  (resp.~$ \, \text{f}_h \mapsto M_{h+1,h} \, \forall\, h=1,\dots,n \, $)
gives an isomorphism among the Lie subalgebra of  $ \gerh $  generated by the
$ \text{e}_h $'s  (resp.~$ \text{f}_h $'s)  and the Lie algebra  $ {\frak n}_+ $
(resp.~$ {\frak n}_+ \, $)  of upper (resp.~lower) triangular square matrix of
size  $ n+1 $;  then for  $ h<k $  the element  $ \text{e}_{h,k} $  corresponds
to the matrix  $ \, {(-1)}^{k-h-1} M_{h,k} \, $,  and it is the root vector
$ \text{e}_\gamma $   --- in the notation of [Ga] ---   associated to the
positive root  $ \, \gamma = \sum_{i=h}^{k-1} \alpha_i \, $  (the
$ \alpha_i $'s  being the simple roots of  $ SL(n+1) $),  and for  $ h>k $  the
element  $ \text{f}_{h,k} $  corresponds to the matrix  $ \, {(-1)}^{k-h-1}
M_{h,k} \, $,  and it is the root vector  $ \text{f}_\gamma $   --- in the
notation of [Ga] ---   associated to the negative root  $ \, -\gamma =
-\sum_{i=k}^{h-1} \alpha_i \, $,  hence corresponding to  $ M_{h,k} $.
Actually, one can also make different choices for such root vectors, but for the
condition that when one of them   --- say  $ \text{e}_\gamma $  ---   is
multiplied by a scalar  $ \, c \in k \setminus \{0\} \, $  then the opposite
one   ---   $ \text{f}_\gamma $  in our case ---   is multiplied by the inverse
scalar  $ c^{-1} \, $;  hence the right-hand-side part in the above formulae
expressing  $ \delta $  does not change.

\vskip1,3truecm

\centerline { \bf  \S \; 2 \,  The quantum function algebra
$ \fqPsln $ }

\vskip10pt

  Let  $ U_q^{\scriptscriptstyle Q} \big({\frak sl}(n+1)\big) $  be the
quantized universal enveloping algebra of Drifel'd and Jimbo (cf.~[Ji] or [DL],
or \S 5.2 later on).  Let  $ \fqPsln $  be its restricted dual Hopf algebra: it
is known (cf.~[APW], Appendix) that  $ \fqPsln $  has the following
presentation:
                              it is\break
 \eject
\noindent   the unital associative  \hbox{$\kq$--alge}bra  generated by
\,  $ \{ \rho_{ij} \mid i, j = 1, \ldots, n+1 \} $  \, with relations
  $$  \eqalignno {
   \rho_{ij} \rho_{ik} = q \, \rho_{ik} \rho_{ij} \; ,  \quad \quad  \rho_{ik}
\rho_{hk}  &  = q \, \rho_{hk} \rho_{ik}  &   \forall\, j<k, i<h  \qquad  \cr
   \rho_{il} \rho_{jk} = \rho_{jk} \rho_{il} \; ,  \quad \quad  \rho_{ik}
\rho_{jl} - \rho_{jl} \rho_{ik}  &  = \left( q - \qm \right) \, \rho_{il}
\rho_{jk}  \hskip90pt  &   \forall\, i<j, k<l  \qquad  \cr
   {det}_q (\rho_{ij}) = 1  &  {}  &   {}  \cr }  $$
where  $ {det}_q $  denotes the so-called  {\it quantum determinant\/},  defined
as
  $$  {det}_q (\rho_{ij}) := \sum_{\sigma \in S_{n+1}}{(-
q)}^{l(\sigma)} \rho_{1,\sigma(1)} \rho_{2,\sigma(2)}
\cdots\rho_{n+1,\sigma(n+1)} \, .  $$
   \indent   The comultiplication  $ \Delta $,  the counit  $ \epsilon $,
and the antipode  $ S $  are given by
  $$  \displaylines{
   \hfill   \Delta (\rho_{ij}) = \sum_{k=1}^n \rho_{ik} \otimes \rho_{kj}
\hfill   \forall\; i, j = 1, \dots, n+1 \, \phantom{.}  \qquad  \cr
   \hfill   \epsilon(\rho_{ij}) = \delta_{ij}   \hfill   \forall\; i, j= 1,
\dots, n+1 \, \phantom{.}  \qquad  \cr
   \hfill   S(\rho_{ij}) = {(-q)}^{j-i} {det}_q \left( {(\rho_{hk})}_{h \neq
j}^{k \neq i} \right)   \hfill   \forall\; i, j = 1, \dots, n+1 \, .  \qquad
\cr }  $$

\vskip1,3truecm

\centerline { \bf  \S \; 3 \,  The integer forms
$ \fqQtildesln $  and  $ \fqPtildesln $ }

\vskip10pt

  {\bf  Definition 3.1.} \  We define  $ \fqQtildesln $  to be the
$ \kqqm $-subalgebra  (with 1) of  $ \fqPsln $  generated by the elements
  $$  \varphi_i := {\, \rho_{ii} - {\rho_{i+1,i+1} \,} \over {\, q - 1 \,}} \; ,
 \qquad  r_{ij} := {\left( q - \qm \right)}^{\delta_{ij}-1} \rho_{ij}   \qquad
\quad  \forall\,  i,j = 1, \dots, n+1.  $$

\vskip7pt

  {\bf 3.2  Presentation of  $ \fqQtildesln $.} \  The presentation of  $
\fqPsln $  above induces a similar presentation of  $ \fqQtildesln $:  it is the
associative  $ \kqqm $--algebra  with 1 given by generators
$ \, \varphi_i \, $,  $ \, r_{ij} \, $,  and relations
  $$  \eqalignno {
%
% relations among the r_{ij} %
%
   r_{ij} r_{ik} = q \, r_{ik} r_{ij} \; ,  \quad \quad  r_{ik} r_{hk}  &  = q
\, r_{hk} r_{ik}  &   \forall\, j<k, i<h  \cr
   r_{il} r_{jk} = r_{jk} r_{il} \; ,  \quad \quad r_{ik} r_{jl} - r_{jl}
r_{ik}  &  = {\left( q - \qm \right)}^{1 + \delta_{ik} + \delta_{jl} -
\delta_{il} - \delta_{jk}} \, r_{il} r_{jk}  \hskip90pt  &   \forall\, i<j,
k<l  \cr
   \widetilde{det}_q (r_{ij})  &  = 1  &   {}  \cr }  $$
(where  $ \widetilde{det}_q $  is defined as
  $$  \widetilde{det}_q \left( {(x_{rs})}_{r,s=1,\dots,N} \right) :=
\sum_{\sigma \in S_N} {(-q)}^{l(\sigma)} {\left( q - \qm \right)}^{e(\sigma)}
x_{1,\sigma(1)} x_{2,\sigma(2)} \cdots x_{N,\sigma(N)}  $$
where  $ \, e(\sigma) := \sum_{t=1}^N (1 - \delta_{t,\sigma(t)} \, $)
  $$  \displaylines{
%
% relation "defining" \varphi_i %
%
   \hfill   (q - 1) \, \varphi_i = r_{ii} - r_{i+1,i+1}  \hskip92pt {}
\hfill   \forall \, i = 1, \dots, n  \cr
%
% commutation relations among the \varphi_i and the r_{jk} %
%
   \hfill   \varphi_i r_{jk} - r_{jk} \varphi_i = 0   \hfill   \forall \,  j<i,
k>i+1 \, ,  \; \forall\, j>i+1, k<i  \cr
   \hfill   \varphi_i r_{jk} - r_{jk} \varphi_i = {( q - 1 )}^{1 + \delta_{jk}}
{\left( 1 + \qm \right)}^{2 + \delta_{jk}} (r_{i+1,k} r_{j,i+1} - r_{ik}
r_{ji})  \quad   \hfill   \forall\, j<i, k<i  \cr }  $$
 \eject
  $$  \eqalignno{
   \varphi_i r_{jk} - r_{jk} \varphi_i = - ( q - 1  &  {)}^{1 + \delta_{jk}}
{\left( 1 + \qm \right)}^{2 + \delta_{jk}} (r_{i+1,k} r_{j,i+1} - r_{ik}
r_{ji})  &   \forall\, j > i \! + \! 1, k > i \! + \! 1  \cr
   \varphi_i r_{ji} - r_{ji} \varphi_i  &  = - r_{ii} r_{ji} + (q-1) {\left( 1 +
\qm \right)}^2 r_{j,i+1} r_{i+1,i}  &   \forall\, j<i  \cr
   {}  &  \varphi_i r_{ji} - r_{ji} \varphi_i = r_{ji} r_{ii}  &   \forall\,
j>i+1  \cr
   \varphi_i r_{j,i+1}  &  - r_{j,i+1} \varphi_i = r_{i+1,i+1} r_{j,i+1}  &
\forall\, j<i  \cr
   \varphi_i r_{j,i+1} - r_{j,i+1} \varphi_i  &  = - r_{j,i+1} r_{i+1,i+1} +
(q-1) {\left( 1 + \qm \right)}^2 r_{i,i+1} r_{ji}  &   \forall\, j<i  \cr
   \varphi_i r_{ij} - r_{ij} \varphi_i  &  = - r_{ii} r_{ij} + (q-1) {\left( 1 +
\qm \right)}^2 r_{i,i+1} r_{i+1,j}  &   \forall\, j<i  \cr
   {}  &  \varphi_i r_{ij} - r_{ij} \varphi_i = r_{ij} r_{ii}  &   \forall\,
j>i+1  \cr
   \varphi_i r_{i+1,j}  &  - r_{i+1,j} \varphi_i = r_{i+1,i+1} r_{i+1,j}  &
\forall\, j<i  \cr
   \varphi_i r_{i+1,j} - r_{i+1,j} \varphi_i =  &  \, - r_{i+1,j} r_{i+1,i+1} +
(q-1) {\left( 1 + \qm \right)}^2 r_{i,j} r_{i+1,i}  &   \forall\, j>i+1  \cr
   \varphi_i r_{ii} -  &  \, r_{ii} \varphi_i = {(q-1)}^2 {\left( 1 + \qm
\right)}^3 r_{i+1,i} r_{i,i+1}  &   \forall\, i  \cr
   \varphi_i r_{i+1,i+1} -  & \, r_{i+1,i+1} \varphi_i = {(q-1)}^2 {\left( 1 +
\qm \right)}^3 r_{i+1,i} r_{i,i+1}  &   \forall\, i  \cr
   \varphi_i r_{i,i+1}  &  - r_{i,i+1} \varphi_i = r_{i,i+1} r_{ii} +
r_{i+1,i+1} r_{i,i+1}  &   \forall\, i  \cr
   \varphi_i r_{i+1,i}  &  - r_{i+1,i} \varphi_i = r_{i+1,i} r_{ii} +
r_{i+1,i+1} r_{i+1,i}  &   \forall\, i  \cr
%
% commutation relations among the \varphi_i %
%
   \varphi_i \varphi_j - \varphi_j \varphi_i =   \hskip51pt  &  {}  &   \cr
   = ( q \! - \! 1) {\left( 1 \! + \! \qm \right)}^3 \big( r_{ij} r_{ji}  &  +
r_{i+1,j+1} r_{j+1,i+1} - r_{i,j+1} r_{j+1,i} - (1 \! - \! \delta_{i+1,j}) \,
r_{i+1,j} r_{j,i+1} \big)   \hskip23pt  &   \forall\, i, j  \cr }  $$
   \indent  Moreover, from the very definitions we also get that
$ \fqQtildesln $  is a Hopf subalgebra of  $ \fqPsln $,  with Hopf structure
uniquely determined by the following formulae:
  $$  \displaylines {
   {} \hfill   \Delta(r_{ij}) = r_{ii} \otimes \, r_{ij} + r_{ij} \otimes r_{jj}
+ (q-1) \left( 1 + \qm \right) \sum_{{k=1} \atop {k \neq i,j}}^{n+1}
r_{ik}\otimes r_{kj}   \hfill  \forall \,  i \neq j  \cr
   {} \hfill   \Delta(r_{ii}) = r_{ii} \otimes r_{ii} + {\left( q - 1 \right)}^2
{\left( 1 + \qm \right)}^2 \sum_{{k=1} \atop {k \neq i}}^{n+1} r_{ik} \otimes
r_{kj}   \hfill  \forall \, i  \cr
   \Delta(\varphi_i) = r_{\scriptscriptstyle i,i} \otimes \varphi_i +
\varphi_i \otimes r_{\scriptscriptstyle i+1,i+1} + (q-1) {\left( 1 + \qm
\right)}^2 \! \left( \sum_{{k=1} \atop {k \neq i}}^{n+1} \!
r_{\scriptscriptstyle i,k} \otimes r_{\scriptscriptstyle k,i} - \!\! \sum_{{k=1}
\atop {k \neq i+1}}^{n+1} \! r_{\scriptscriptstyle i+1,k} \otimes
r_{\scriptscriptstyle k,i+1} \! \right)   \hfill  \forall \, i  \cr
   {} \hfill   S (r_{ij}) = {(-q)}^{j-i} \widetilde{det}_q \left( {(r_{hk})}_{h
\neq j}^{k \neq i} \right)   \hfill  \forall \, i, j  \cr
   {} \quad  S(\varphi_i) = - r_{1,1} r_{2,2} \cdots r_{i-1,i-1} \varphi_i
r_{i+2,i+2} \cdots r_{n+1,n+1} +   \hfill   \cr
   \hfill    + \sum_{\sigma \in S_n \setminus \{1\}} {(-q)}^{l(\sigma)} {\left(
q - 1 \right)}^{e(\sigma) - 1} {\left( 1 + \qm \right)}^{e(\sigma)} \left(
\prod_{{j=1} \atop {j \neq i+1}}^{n+1} r_{j,\sigma(j)} - \prod_{{j=1} \atop {j
\neq i}}^{n+1} r_{j,\sigma(j)} \right)   \hfill  \forall \,  i  \cr
   \hfill  \varepsilon (r_{ij}) = \delta_{ij} \, ,  \qquad  \varepsilon
(\varphi_i) = 0  \hfill   \forall \,  i, j  \cr }  $$

\vskip7pt

  {\bf Remark 3.3.} \  It is clear by definition that
$ \fqQtildesln $  is a  $ \kqqm $--integer  form of  $ \fqPsln $:  in other
words, it is a Hopf  $ \kqqm $--subalgebra  of  $ \fqPsln $  which is flat as
a  $ \kqqm $--module and is such that  $ \, \kq \otimes_{k \, \left[ q, \qm
\right]} \fqQtildesln \cong \fqPsln \, $  as Hopf  $ \kq $--algebras.

\vskip7pt

  {\bf  Definition 3.4.} \  We define  $ \fqPtildesln $  to be the
$ \kqqm $-subalgebra  (with 1) of  $ \fqPsln $  generated by the elements
  $$  \psi_i := {{\, \rho_{1{}1} \rho_{2{}2} \cdots \rho_{ii} - 1 \,} \over
{\, q - 1 \,}} \; ,  \qquad  r_{ij} := {\left( q - \qm
\right)}^{\delta_{ij}-1} \rho_{ij}   \qquad \quad  \forall\,  i,j = 1,
\dots, n+1.  $$

\vskip7pt

  {\bf 3.5  Presentation of  $ \fqPtildesln $.} \  The presentation of  $
\fqPsln $  above induces a similar presentation of  $ \fqPtildesln $:  in fact
the latter is the associative  $ \kqqm $--algebra  with 1 given by generators
$ \, \psi_i \, $,  $ \, r_{ij} \, $,  and relations
  $$  \displaylines {
%
% relations among the r_{ij} %
%
   \hfill   r_{ij} r_{ik} = q \, r_{ik} r_{ij} \; ,  \qquad \quad  r_{ik} r_{hk}
= q \, r_{hk} r_{ik}   \hfill  \forall\, j<k, i<h  \;  \cr
   \hfill   r_{il} r_{jk} = r_{jk} r_{il} \; ,  \quad \quad r_{ik} r_{jl} -
r_{jl} r_{ik} = {\left( q - \qm \right)}^{1 + \delta_{ik} + \delta_{jl} -
\delta_{il} - \delta_{jk}} \, r_{il} r_{jk}   \hfill  \forall\, i<j, k<l \;  \cr
   \widetilde{det}_q (r_{ij}) = 1  \cr
%
% relation "defining" \psi_i %
%
   \hfill   (q-1) \, \psi_i = r_{1{}1} r_{2{}2} \cdots r_{ii} - 1   \hfill
\forall \, i= 1, \dots, n+1 \;  \cr
%
% commutation relations among the \psi_i and the r_{jk} %
%
   {} \qquad   \psi_i r_{jk} = q^{1 - \eta(i,j,k) - \zeta(i,j,k)} r_{jk} \psi_i
+ \theta(i,j,k) \cdot r_{jk} + {\left( q - 1 \right)}^{1+\delta_{jk}}
{\left( 1 + \qm \right)}^{2+\delta_{jk}} \cdot  \hfill {\ }  \cr
   {\ \ }   \hfill  \cdot \sum_{s=1}^i \big( \eta(i,j,k) - \zeta(i,j,k)
\big) r_{1{}1} \cdots r_{s-1,s-1} r_{sk} r_{js} r_{s+1,s+1} \cdots r_{i{}i}
\hfill   \forall \, i, j, k \;  \cr }  $$
(where  $ \quad \theta(i,j,k):= 1 \; \forall \, j \leq i < k \, $  or
$ \, k \leq i < j \, $,  $ \; \eta(i,j,k):= 1 \; \forall \, i < j \wedge k
\, $,  $ \; \zeta(i,j,k):= 1 $\break
$ \forall \, i \geq j \vee k \, $,  \; whilst
$ \, \theta(i,j,k) $,  $ \, \eta(i,j,k) \, $  and  $ \, \zeta(i,j,k) \, $
are zero in the other cases)
%
% special relation for \psi_{n+1}
%
  $$  \psi_{n+1} =
- \sum_{ \Sb
           \sigma \in S_{n+1}  \\
           \sigma \neq id  \\
         \endSb } {(-q)}^{l(\sigma)} {\left( q - 1 \right)}^{e(\sigma)-1}
{\left( 1 + \qm \right)}^{e(\sigma)} r_{1,\sigma(1)} r_{2,\sigma(2)} \cdots
r_{n+1,\sigma(n+1)}  $$
%
% commutation relations among the \psi_i %
%
  $$  \displaylines {
   {} \quad  \psi_i \psi_j - \psi_j \psi_i = (q-1) {\left( 1 + \qm \right)}^3
\cdot  \hfill   {\ }  \cr
   {\ } \hfill   \cdot \sum_{\scriptscriptstyle k=i+1}^{\scriptscriptstyle j} \!
\sum_{\scriptscriptstyle s=1}^{\scriptscriptstyle i} r_{\scriptscriptstyle 1{}1}
r_{\scriptscriptstyle 2{}2} \cdots r_{\scriptscriptstyle k-1, k-1} \cdot
r_{\scriptscriptstyle 1{}1} r_{\scriptscriptstyle 2{}2} \cdots
r_{\scriptscriptstyle s-1,s-1} r_{\scriptscriptstyle sk} r_{\scriptscriptstyle
ks} r_{\scriptscriptstyle s+1,s+1} \cdots r_{\scriptscriptstyle ii} \cdot
r_{\scriptscriptstyle k+1,k+1} r_{\scriptscriptstyle k+2,k+2} \cdots
r_{\scriptscriptstyle jj}  \, \forall\, i \! < \! j  \cr }  $$
   \indent   Furthermore, from the very definitions we also get that
$ \fqPtildesln $  is a Hopf subalgebra of  $ \fqPsln $,  with Hopf structure
uniquely determined by the following formulae:
  $$  \displaylines {
   \hfill   \Delta(r_{ij}) = r_{ii} \otimes r_{ij} +
r_{ij} \otimes r_{jj} + (q-1) \left( 1 + \qm \right) \sum_{{k=1} \atop {k \neq
i,j}}^{n+1} r_{ik}\otimes r_{kj}   \hfill  \forall \,  i \neq j \;  \cr
   \hfill   \Delta(r_{ii}) = r_{ii} \otimes r_{ii} + {\left( q - 1 \right)}^2
{\left( 1 + \qm \right)}^2 \sum_{{k=1} \atop {k \neq i}}^{n+1} r_{ik} \otimes
r_{kj}   \hfill  \forall \, i \;  \cr }  $$
  $$  \Delta(\psi_i) \! = \! \left( 1 + \qm \right) \cdot \! \sum_s \! {\left( q
- \qm \right)}^{2 \cdot N(s) - 1} \! \prod_{k=1}^i r_{\scriptscriptstyle k,s(k)}
\! \otimes \! r_{\scriptscriptstyle s(k),k} + \psi_{\scriptscriptstyle i} \!
\otimes \! r_{\scriptscriptstyle 1,1} r_{\scriptscriptstyle 2,2} \cdots
r_{\scriptscriptstyle i,i} + 1 \! \otimes \! \psi_{\scriptscriptstyle i}
\eqno  \forall \, i \;  $$
(where  $ s $  ranges over all maps  $ \, s \colon \{1,2,\dots,i\} \rightarrow
\{1,2,\dots,n+1\} \, $  such that  $ \, s(j) \neq j \, $  for some  $ j \in
\{1,2,\dots,i\} $,  and  $ \, N(s) := \sum_{j=1}^i \big( 1 - \delta_{j,s(j)}
\big) \, $)
  $$  \displaylines {
   \hfill   S (r_{ij}) = {(-q)}^{j-i} \widetilde{det}_q \left( {(r_{hk})}_{h
\neq j}^{k \neq i} \right)  \hfill   \forall \,  i, j  \cr
   \hfill   S(\psi_i) = - \psi_i + {\Cal O}(q-1)   \hfill  \forall \, i  \cr
   \hfill  \varepsilon (r_{ij}) = \delta_{ij} \, ,  \qquad  \varepsilon (\psi_i)
= 0  \hfill   \forall \,  i, j  \cr }  $$
where  $ {\Cal O}(q-1) $  denotes some element of  $ \, (q-1) \cdot \fqPtildesln
\, $.  To give an example, we show that  $ S(\psi_i) = - \psi_i +
{\Cal O}(q-1) $.  By definition we have
  $$  S(\psi_i) = S \left( { \, r_{1{}1} r_{2{}2} \cdots r_{ii} - 1 \, \over \,
q-1 \, } \right) = { \, S \left( r_{i{}i} \right) \cdots S \left( r_{2{}2}
\right) S \left( r_{1{}1} \right) - 1 \, \over \, q-1 \, } \; ;  $$
but
  $$  \displaylines{
   {\ }  S \left( r_{j{}j} \right) = \widetilde{det}_q \left( {(r_{hk})}_{h,k
\neq j} \right) =   \hfill {\ }  \cr
   {} \hfill   = r_{1{}1} r_{2{}2} \cdots r_{j-1,j-1} \cdot r_{j+1,j+1} \cdots
r_{n+1,n+1} + {\Cal O} \left( {\left( q-1 \right)}^2 \right) = \prod_{\Sb  s = 1
\\  s \neq j \\  \endSb}^{n+1} r_{s{}s} + {\Cal O} \left( {\left( q-1 \right)}^2
\right)  \cr }  $$
for all  $ j $,  and
  $$  1 = \widetilde{det}_q \left( r_{hk} \right) = r_{1{}1} r_{2{}2} \cdots
r_{n+1,n+1} + {\Cal O} \left( {\left( q-1 \right)}^2 \right)  $$
therefore
  $$  S (\psi_i) = { \, \prod_{\Sb  s = 1 \\  s \neq i \\  \endSb}^{n+1}
r_{s{}s} \cdot \prod_{\Sb  s = 1 \\  s \neq i-1 \\  \endSb}^{n+1} r_{s{}s}
\cdots \prod_{\Sb  s = 1 \\  s \neq 2 \\  \endSb}^{n+1} r_{s{}s} \cdot
\prod_{\Sb  s = 1 \\  s \neq 1 \\  \endSb}^{n+1} r_{s{}s} - {\left(
\prod_{s=1}^{n+1} r_{s{}s} \right)}^i \, \over \, q-1 \, } + {\Cal O} \left( q-1
\right) \; ;  $$
now using the fact that  $ \, r_{h{}h} r_{k{}k} = r_{k{}k} r_{h{}h} + {\Cal O}
\left( {\left( q-1 \right)}^3 \right) \, $  we get
  $$  S (\psi_i) = \prod_{\Sb  s = 1 \\  s \neq i \\  \endSb}^{n+1} r_{s{}s}
\cdot \prod_{\Sb  s = 1 \\  s \neq i-1 \\  \endSb}^{n+1} r_{s{}s} \cdots
\prod_{\Sb  s = 1 \\  s \neq 2 \\  \endSb}^{n+1} r_{s{}s} \cdot
\prod_{\Sb  s = 1 \\  s \neq 1 \\  \endSb}^{n+1} r_{s{}s} \cdot { \, 1 -
r_{1{}1} r_{2{}2} \cdots r_{i{}i} \, \over \, q-1 \, } + {\Cal O} \left( q-1
\right)  $$
and finally, since  $ \, r_{j{}j} = 1 + {\Cal O} \left( q-1 \right) \, $  (as
one easily gets from relations  $ \, (q-1) \, \psi_s = r_{1{}1} r_{2{}2} \cdots
r_{s{}s} - 1 \, $),  we find
  $$  S (\psi_i) = { \, 1 - r_{1{}1} r_{2{}2} \cdots r_{i{}i} \, \over \, q-1 \,
} + {\Cal O} \left( q-1 \right) = -\psi_i + {\Cal O} \left( q-1 \right) \, ,
\quad \hbox{q.e.d.} $$

\vskip7pt

  {\bf Remark 3.6.} \  Here again, it is clear   --- by definition and by the
description of the Hopf structure ---   that  $ \fqPtildesln $  is a
$ \kqqm $--integer  form of  $ \fqPsln $.

\vskip7pt

  {\bf  Definition 3.7.} \  We define  $ \fqtildesln $  to be the
$ \kqqm $-subalgebra  (with 1) of  $ \fqPsln $  generated by the elements
  $$  \chi_i := {{\, \rho_{i{}i} - 1 \,} \over {\, q - 1 \,}} \; ,  \qquad
r_{ij} := {\left( q - \qm \right)}^{\delta_{ij}-1} \rho_{ij}   \qquad \quad
\forall\,  i,j = 1, \dots, n+1.  $$

\vskip7pt

  {\bf 3.8  Presentation of  $ \fqtildesln $.} \  Again, we have a presentation
of  $ \fqtildesln $:  it is the associative  $ \kqqm $--algebra  with 1 given by
generators  $ \, \chi_i \, $,  $ \, r_{ij} \, $,  and relations
  $$  \displaylines {
%
% relations among the r_{ij} %
%
   \hfill   r_{ij} r_{ik} = q \, r_{ik} r_{ij} \; ,  \qquad \quad  r_{ik} r_{hk}
= q \, r_{hk} r_{ik}   \hfill  \forall\, j<k, i<h \;  \cr
   \hfill   r_{il} r_{jk} = r_{jk} r_{il} \; ,  \qquad  r_{ik} r_{jl} - r_{jl}
r_{ik} = {\left( q - \qm \right)}^{1 + \delta_{ik} + \delta_{jl} - \delta_{il} -
\delta_{jk}} \, r_{il} r_{jk}   \hfill  \forall\, i<j, k<l \;  \cr
   \widetilde{det}_q (r_{ij}) = 1  \hskip89pt  \cr
%
% relation "defining" \chi_i %
%
   \hfill   (q-1) \, \chi_i = r_{ii} - 1   \hfill  \forall \,  i= 1, \dots, n+1
\;  \cr
%
% commutation relations among the \chi_i and the r_{jk} %
%
   \hfill   \qquad  \chi_i r_{jk} - r_{jk} \chi_i = 0   \hfill  \forall \, j<i<k
\, ,  \; \forall\, j>i>k \;  \cr
   \hfill   \chi_i r_{ji} - r_{ji} \chi_i = - r_{ii} r_{ji}   \hfill  \qquad
\qquad \forall \, j<i \;  \cr
   \hfill   \chi_i r_{ji} - r_{ji} \chi_i = + r_{ii} r_{ji}   \hfill  \qquad
\qquad \forall \, j>i \;  \cr
   \hfill   \chi_i r_{ii} - r_{ii} \chi_i = 0   \hfill  \quad \qquad \qquad
\qquad \forall \, i \;  \cr
   \hfill   \chi_i r_{ik} - r_{ik} \chi_i = - r_{ii} r_{ik}   \hfill  \qquad
\qquad  \forall \, k<i \;  \cr
   \hfill   \chi_i r_{ik} - r_{ik} \chi_i = + r_{ii} r_{ik}   \hfill  \qquad
\qquad  \forall \, k>i \;  \cr
   \hfill   \chi_i r_{jk} - r_{jk} \chi_i = - {(q-1)}^2 {\left( 1 + \qm
\right)}^3 r_{ii} r_{ik}   \hfill  \forall \, j,k<i \;  \cr
   \hfill   \chi_i r_{jk} - r_{jk} \chi_i = + {(q-1)}^2 {\left( 1 + \qm
\right)}^3 r_{ii} r_{ik}   \hfill  \forall \, j,k>i \;  \cr
%
% commutation relations among the \chi_i %
%
   \hfill   \chi_i \chi_j - \chi_j \chi_i = (1 - \delta_{ij}) \cdot (q-1)
{\left( 1 + \qm \right)}^3 r_{ij} r_{ji}   \hfill  \forall \, i \leq j \;  \cr
%
% special relation for the \chi_i %
%
   {\ } \quad   \sum_{i=1}^{n+1} r_{1,1} r_{2,2} \cdots r_{i-1,i-1} \chi_i =
\hfill {\ }  \cr
   {\ } \hfill   = \sum_{\sigma \in S_{n+1} \setminus
\{1\}} {(-q)^{l(\sigma)}} \, {(q-1)}^{e(\sigma)-1} {\left( 1 + \qm
\right)}^{e(\sigma)} r_{1,\sigma(1)} r_{2,\sigma(2)} \cdots r_{n+1,\sigma(n+1)}
 \quad {}  \cr }  $$
   \indent  A straightforward verification shows that
$ \fqtildesln $  is also a Hopf subalgebra of  $ \fqPsln $,  whose Hopf
structure is given by formulae
  $$  \displaylines {
   \hfill   \Delta(r_{ij}) = r_{ii} \otimes r_{ij} + r_{ij} \otimes r_{jj} +
(q-1) \left( 1 + \qm \right) \sum_{{k=1} \atop {k \neq i,j}}^{n+1} r_{ik}
\otimes r_{kj}   \hfill  \forall \,  i \neq j \;  \cr
   \hfill   \Delta(r_{ii}) = r_{ii} \otimes r_{ii} + {\left( q - 1 \right)}^2
{\left( 1 + \qm \right)}^2 \sum_{{k=1} \atop {k \neq i}}^{n+1} r_{ik} \otimes
r_{kj}   \hfill  \forall \, i \;  \cr }  $$
  $$  \displaylines{
   \hfill   \Delta(\chi_i) = r_{ii} \otimes \chi_i + (q-1) {\left( 1 + \qm
\right)}^2 \sum_{{k=1} \atop {k \neq i}}^{n+1} r_{i,k} \otimes r_{k,i}
\hfill  \forall \, i \;  \cr
   \hfill   S (r_{ij}) = {(-q)}^{j-i} \widetilde{det}_q \left( {(r_{hk})}_{h
\neq j}^{k \neq i} \right)   \hfill  \forall \, i, j \;  \cr
   {} \quad  S(\chi_i) = - r_{1,1} r_{2,2} \cdots r_{i-1,i-1} \chi_i r_{i+2,i+2}
\cdots r_{n+1,n+1} +   \hfill   \cr
   + \sum_{\sigma \in S_n \setminus \{1\}} {(-q)}^{l(\sigma)} {\left( q - 1
\right)}^{e(\sigma) - 1} {\left( 1 + \qm \right)}^{e(\sigma)} \prod_{\Sb  {j=1}
\\  j \neq i \\  \sigma(j) \neq i \\  \endSb}^{n+1} r_{j,\sigma(j)} -  \cr
   {\ } \hfill   - \sum_{\sigma \in S_{n+1} \setminus \{1\}} {(-q)}^{l(\sigma)}
{\left( q - 1 \right)}^{e(\sigma) - 1} {\left( 1 + \qm \right)}^{e(\sigma)}
\prod_{j=1}^{n+1} r_{j,\sigma(j)}   \qquad  \forall \, i \;  \cr
   \hfill  \varepsilon (r_{ij}) = \delta_{ij} \, ,  \qquad \qquad  \varepsilon
(\chi_i) = 0   \qquad  \hfill   \forall  \,  i, j \, . \;  \cr }  $$

\vskip7pt

  {\bf Remark 3.9.} \  By definition and by the description of its Hopf
structure we see at once that  $ \fqtildesln $  is a  $ \kqqm $--integer  form
of  $ \fqPsln $.

\vskip1,3truecm

\centerline { \bf  \S \; 4 \,  The main theorem: specialization results }

\vskip10pt

   Since the integer forms we are dealing with are Hopf algebras  {\sl over the
ring  $ \kqqm $},  we can consider their specialization at  $ \, q=1 \, $,
namely the Hopf  $ k $--algebras
  $$  \eqalign{
   \funoQtildesln  &  := \fqQtildesln \Big/ (q-1) \, \fqQtildesln \, ,  \cr
   \funoPtildesln  &  := \fqPtildesln \Big/ (q-1) \, \fqPtildesln \, ,  \cr
    \funotildesln  &  := \fqtildesln \Big/ (q-1) \, \fqtildesln \, .  \cr }  $$
   \indent   Our main result is the following:

\proclaim{Theorem 4.1} \  The Hopf  $ k $--algebras  $ \funoQtildesln $,
$ \funoPtildesln $,  and  $ \funotildesln $,  are Poisson Hopf coalgebras
isomorphic to  $ \uh $.  In other words,  $ \fqQtildesln $,  $ \fqPtildesln $,
and  $ \fqtildesln $  all specialize to the Poisson Hopf coalgebra  $ \uh $.
\endproclaim

\demo{Proof}  Consider  $ \funoQtildesln $;  it inherits from  $ \fqQtildesln $
the following presentation (which is obtained from that of  $ \fqQtildesln $  by
setting  $ \, q = 1 \, $):  it is the unital associative  $ k $--algebra  with
generators  $ r_{ij} $,  $ \varphi_k $  ($ \, i, j = 1, \ldots, n+1 $;
$ k = 1, \dots, n \, $)  and relations
  $$  \displaylines {
   \hfill   r_{ij} r_{ik} = r_{ik} r_{ij} \; ,  \qquad \quad  r_{ik} r_{hk} =
r_{hk} r_{ik}   \hfill  \forall\, j<k, i<h \;  \cr
   \hfill   r_{il} r_{jk} = r_{jk} r_{il} \; ,  \quad \quad  r_{ik} r_{jl} -
r_{jl} r_{ik} = {(0)}^{1 + \delta_{ik} + \delta_{jl} - \delta_{il} -
\delta_{jk}} \, r_{il} r_{jk}   \hfill  \forall\, i<j, k<l \;  \cr
   r_{1,1} r_{2,2} \cdots r_{n+1,n+1} = 1  \hskip60pt  \cr
   \hfill   0 = r_{ii} - r_{i+1,i+1}   \hfill  \forall \, i = 1, \dots, n \;
\cr }  $$
  $$  \displaylines{
   \hfill   \qquad \qquad \qquad  \varphi_i r_{jk} - r_{jk} \varphi_i = 0
\hfill  \forall \, j<i, k>i+1 \, ,  \; \forall\, j>i+1, k<i \;  \cr
   \hfill   \varphi_i r_{jk} - r_{jk} \varphi_i = 0  \quad \qquad   \hfill
\forall \, j<i, k<i \;  \cr
   \hfill   \varphi_i r_{jk} - r_{jk} \varphi_i = 0   \hfill  \forall \, j>i+1,
k>i+1 \;  \cr
   \hfill   \varphi_i r_{ji} - r_{ji} \varphi_i = - r_{ii} r_{ji}  \qquad \qquad
  \hfill  \forall \, j<i \;  \cr
   \hfill   \varphi_i r_{ji} - r_{ji} \varphi_i = r_{ji} r_{ii}  \qquad \quad
\hfill  \forall \, j>i+1 \;  \cr
   \hfill   \varphi_i r_{j,i+1} - r_{j,i+1} \varphi_i = r_{i+1,i+1} r_{j,i+1}
\quad   \hfill  \forall \, j<i \;  \cr
   \hfill   \varphi_i r_{j,i+1} - r_{j,i+1} \varphi_i = - r_{j,i+1}
r_{i+1,i+1}   \hfill  \forall \, j<i \;  \cr
   \hfill   \varphi_i r_{ij} - r_{ij} \varphi_i = - r_{ii} r_{ij}  \quad
\qquad   \hfill  \forall \, j<i \;  \cr
   \hfill   \varphi_i r_{ij} - r_{ij} \varphi_i = r_{ij} r_{ii}  \qquad   \hfill
 \forall \, j>i+1 \;  \cr
   \hfill   \varphi_i r_{i+1,j} - r_{i+1,j} \varphi_i = r_{i+1,i+1} r_{i+1,j}
\quad   \hfill  \forall \, j<i \;  \cr
   \hfill   \; \varphi_i r_{i+1,j} - r_{i+1,j} \varphi_i = - r_{i+1,j}
r_{i+1,i+1}   \hfill  \forall \, j>i+1 \;  \cr
   \hfill   \varphi_i r_{ii} - r_{ii} \varphi_i = 0  \qquad \qquad \quad
\hfill  \forall \, i \;  \cr
   \hfill   \varphi_i r_{i+1,i+1} - r_{i+1,i+1} \varphi_i = 0  \qquad \qquad
\hfill  \forall \, i \;  \cr
   \hfill   \varphi_i r_{i,i+1} - r_{i,i+1} \varphi_i = r_{i,i+1} r_{ii} +
r_{i+1,i+1} r_{i,i+1}  \qquad   \hfill  \forall \, i \;  \cr
   \hfill   \varphi_i r_{i+1,i} - r_{i+1,i} \varphi_i = r_{i+1,i} r_{ii} +
r_{i+1,i+1} r_{i+1,i}  \qquad   \hfill  \forall \, i \;  \cr
   \hfill   \varphi_i \varphi_j - \varphi_j \varphi_i = 0  \quad \qquad   \hfill
 \forall \, i, j \;  \cr }  $$
   \indent  Moreover, its Hopf structure is given by
  $$  \displaylines {
   \hfill   \Delta(r_{ij}) = r_{ii} \otimes r_{ij} + r_{ij} \otimes r_{jj}
\hfill  \forall \, i \neq j \;  \cr
   \hfill   \Delta(r_{ii}) = r_{ii} \otimes r_{ii}   \hfill  \forall \, i \;
\cr
   \hfill   \Delta(\varphi_i) = r_{ii} \otimes \varphi_i + \varphi_i \otimes
r_{i+1,i+1}   \hfill  \forall \, i \;  \cr
   \hfill   S (r_{ij}) = - r_{ij} \cdot \prod_{k=1 \atop k \neq i,j}^{n+1}
r_{k,k}   \hfill  \forall \, i, j \;  \cr
   \hfill   S(\varphi_i) = - r_{1,1} r_{2,2} \cdots r_{i-1,i-1} \varphi_i
r_{i+2,i+2} \cdots r_{n+1,n+1}   \hfill  \forall \, i \;  \cr
   \hfill  \varepsilon (r_{ij}) = \delta_{ij} \, ,  \qquad  \varepsilon
(\varphi_i) = 0  \hfill   \forall \,  i, j \, . \;  \cr }  $$
   \indent   In particular we have that  $ \, r_{1,1} = r_{2,2} = \cdots =
r_{n+1,n+1} \, $,  hence  $ \, r_{i,i}^{\,n+1} = 1 \, $  for all $ i $,  whence
$ \, r_{i,i} \in k \, $  ($ k $  is algebraically closed); but then  $ \,
\Delta(r_{ii}) = r_{ii} \otimes r_{ii} \, $  implies  $ \, r_{ii} = 1 \, $.
Now relations  $ \, r_{ik} r_{jl} - r_{jl} r_{ik} = {(0)}^{1 + \delta_{ik} +
\delta_{jl} - \delta_{il} - \delta_{jk}} \, r_{il} r_{jk} \, $  ($ i<j, k<l $)
gives in particular  $ \, r_{i,j-1} r_{j-1,j} - r_{j-1,j} r_{i,j-1} =
{(0)}^{\delta_{i,j-1}} \, r_{i,j} r_{j-1,j-1} \, $  ($ i<j $),  and similarly
$ \, r_{i,j+1} r_{j+1,j} - r_{j+1,j} r_{i,j+1} = {(0)}^{\delta_{i,j+1}} \,
r_{i,j} r_{j+1,j+1} \, $  ($ i<j $),  whence we deduce that the elements
$ r_{i,i+1} $,  $ r_{i+1,i} $  ($ i=1, \dots, n $)  together with the
$ \varphi_j $'s  are enough to generate  $ \funoQtildesln $.
                                                  \par
   Now from the relations above one finds that for the generators
$ r_{i,i+1} $,  $ \varphi_i $,  $ r_{i+1,i} $  ($ i=1, \dots, n $)  exactly the
same relations hold than we have in \S 1 for the generators  $ -\text{f}_i $,
$ \text{h}_i $,  $ \text{e}_i $,  of  $ \uh $:  therefore
$ \funoQtildesln $  and  $ \uh $,  having the same presentation, are
isomorphic as  $ k $--algebras.  Nevertheless, the formulae for the values of
Hopf operations (of  $ \funoQtildesln $)  on the generators  $ r_{i,i+1} $,
$ \varphi_i $,  $ r_{i+1,i} $  are exactly the same   --- when taking into
account that  $ \, r_{ii} = 1 \, $  for all  $ i $ ---   than similar formulae
for the generators  $ -\text{f}_i $,  $ \text{h}_i $,  $ \text{e}_i $,  of
$ \uh $: thus the  $ k $--algebra  isomorphism  $ \; \Phi : \funoQtildesln
\longrightarrow \uh \; $  given by
  $$  \Phi : \; \quad  r_{i,i+1} \mapsto -\text{f}_i \, ,  \quad  \varphi_i
\mapsto \text{h}_i \, ,  \quad  r_{i+1,i} \mapsto +\text{e}_i \, ,
\eqno (4.1)  $$
is even one of Hopf algebras.  In particular, the comultiplication of
$ \funoQtildesln $  is cocommutative: hence a Poisson cobracket  $ \, \delta
\colon \funoQtildesln \longrightarrow \funoQtildesln \otimes \funoQtildesln
\, $  is canonically defined by  $ \, \delta \left( x{\big\vert}_{q=1} \right)
:= {\; \left( \Delta - {\Delta}^{op} \right)(x) \; \over \; q-1
\;}{\bigg\vert}_{q=1} \, $.
                                                  \par
   In order to compare the latter Poisson cobracket with the one on  $ \uh $
given from scratch, we have to unravel the preimage in  $ \funoQtildesln $  of
root vectors in  $ \uh $.
                                                \par
  We already saw that  $ \, r_{i{}i} = 1 \, $  ($ \, i= 1, \dots, n+1 \, $)  in
$ \funoQtildesln $;  from this fact and the relations  $ \, r_{ik} r_{jl} -
r_{jl} r_{ik} = {(0)}^{1 + \delta_{ik} + \delta_{jl} - \delta_{il} -
\delta_{jk}} \, r_{il} r_{jk} \, $  ($ \, i < j $,  $ k < l \, $),  for  $ \, j
= k \, $  one gets
  $$  r_{il} = + \big[ r_{ij} \, , \, r_{jl} \big]   \qquad \quad \forall \; i <
j < l \, , $$
and similarly for  $ \, i = l \, $
  $$  r_{jk} = - \big[ r_{ji} \, , \, r_{ik} \big]  \qquad \quad \forall \; j >
i > k \, ;  $$
in particular
  $$  \hbox{ $ \eqalign{
   r_{ij}  &  = + \big[ r_{i,j-1} \, , \, r_{j-1,j} \big] = - \big[ r_{j-1,j} \,
, r_{i,j-1} \big]  \qquad \qquad  \forall\, i < j-1 \, ,  \cr
   r_{ji}  &  = - \big[ r_{j,j-1} \, , \, r_{j-1,i} \big] = + \big[ r_{j-1,i} \,
, \, r_{j,j-1} \big]  \qquad \qquad  \forall\, j > i+1 \, .  \cr } $ }
\eqno (4.2)  $$
   \indent   Comparing (4.2) and (1.1), by a simple induction one gets from
(4.1) that
  $$  \Phi \big( r_{ij} \big) = {(-1)}^{j-i} \, \text{f}_{ji} \; ,  \; \qquad
\Phi \big( r_{ji} \big) = {(-1)}^{j-i-1} \text{e}_{ij} \; ,  \qquad \quad
\forall \; i < j \, .   \eqno (4.3)  $$
   \indent   Now, a straightforward computation gives
  $$  \displaylines{
   \delta(r_{i,i+1}) = \varphi_i \otimes r_{i,i+1} - r_{i,i+1} \otimes \varphi_i
+ 2 \cdot {\sum_{j=1}^{n+1}{}^{\widehat{i}, \, \widehat{i+1}}} \big( r_{i,j}
\otimes r_{j,i+1} - r_{j,i+1} \otimes r_{i,j} \big)  \cr
   \delta(\varphi_i) = 4 \cdot \Bigg( \sum_{j=1}^{n+1} \big( r_{i,j} \otimes
r_{j,i} - r_{j,i} \otimes r_{i,j} \big) - \sum_{j=1}^{n+1} \big( r_{i+1,j}
\otimes r_{j,i+1} - r_{j,i+1} \otimes r_{i+1,j} \big) \Bigg)  \cr
   \delta(r_{i+1,i}) = r_{i+1,i} \otimes \varphi_i - \varphi_i \otimes r_{i+1,i}
+ 2 \cdot {\sum_{j=1}^{n+1}{}^{\widehat{i}, \, \widehat{i+1}}} \big( r_{i+1,j}
\otimes r_{j,i} - r_{j,i} \otimes r_{i+1,j} \big)  \cr }  $$
(for all  $ \, i = 1, \dots, n \, $),  where a superscript  $ \widehat{h} $
means that the index  $ h $  must be discarded; then it is a simple task of
rewriting (using (4.3)) to see that these formulae correspond   --- via
$ \Phi $  ---   to the analogous ones for  $ \uh $.  Thus the isomorphism  $ \;
\Phi : \funoQtildesln \longrightarrow \uh \; $  above is one of Poisson Hopf
coalgebras; so we have proved the claim for  $ \funoQtildesln $.
                                                \par
   As for the other two algebras, we shall shortly conclude relying on the first
one.  In fact  $ \funoPtildesln $  and  $ \funotildesln $  differ from
$ \funoQtildesln $  only for "toral" generators: the  $ \psi_i $'s,  resp.~the
$ \chi_i $'s,  instead  of the  $ \varphi_i $'s.  Now, definitions give at once
  $$  \psi_i = \chi_1 + r_{1,1} \chi_2 + r_{1,1} r_{2,2} \chi_3 + \cdots +
r_{1,1} r_{2,2} \cdots r_{i-1,i-1} \chi_i  \qquad   \eqno \forall \; i=1, \dots,
n+1  \qquad  \;  $$
so that, since  $ \, r_{j,j} \equiv 1 \; \mod\, (q-1) \, $,  we have that  $ \,
\{ \psi_1, \dots, \psi_{n+1} \} $  modulo  $ (q-1) \, $  and  $ \, \{ \chi_1,
\dots, \chi_{n+1} \} $  modulo  $ (q-1) \, $  span the same  $ k $--vector
space, whence
  $$  \funoPtildesln = \funotildesln \, .  $$
   \indent   Furthermore, definitions give also
  $$  \qquad  \varphi_i = \chi_i - \chi_{i+1}  \qquad  \qquad  \forall \; i=1,
\dots, n  $$
hence the  $ k $--span  of  $ \, \{ \varphi_1, \dots, \varphi_{n+1} \} $
modulo  $ (q-1) \, $  is contained in the  $ k $--span of  $ \, \{ \chi_1,
\dots,$  $ \chi_{n+1} \} $  modulo  $ (q-1) \, $;  moreover, the relation
  $$  \displaylines{
   {} \quad   \sum_{i=1}^{n+1} r_{1,1} r_{2,2} \cdots r_{i-1,i-1} \chi_i =
\hfill {\ }  \cr
   {\ } \hfill   = \sum_{\sigma \in S_{n+1} \setminus \{1\}} {(-q)^{l(\sigma)}}
{(q-1)}^{e(\sigma)-1} {\left( 1 + \qm \right)}^{e(\sigma)} r_{1,\sigma(1)}
r_{2,\sigma(2)} \cdots r_{n+1,\sigma(n+1)}   \quad {}  \cr }  $$
for  $ \, q=1 \, $  turns into
  $$  \chi_1 + \cdots + \chi_{n+1} = 0 \; ;  $$
thus the  $ k $--span  of  $ \, \{ \varphi_1, \dots, \varphi_{n+1} \} $  modulo
$ (q-1) \, $  and the  $ k $--span of  $ \, \{ \chi_1, \dots, \chi_{n+1} \} $
modulo  $ (q-1) \, $  have both dimension  $ n $,  hence they coincide.  We
conclude that
  $$  \funoQtildesln = \funotildesln = \funoPtildesln  $$
whence the claim.   $ \square $
\enddemo

\vskip1,3truecm

\centerline { \bf  \S \; 5 \,  $ \widetilde{F}_q^{\scriptscriptstyle M}
[SL(n+1)] $  as approximation of  $ {\Cal F}^{\scriptscriptstyle M}
[SL(n+1)] $ }

\vskip10pt

  {\bf 5.1  Motivations.} \  To explain the definitions of the integer forms of
\S 4 some comments are in order.  We resume the analysis in [Ga], using the same
notation, and make it more explicit for  $ \, G = SL(n+1) \, $.
                                                       \par
   Given the quantized universal enveloping algebra
$ \uqPsln $,  resp.~$ \uqQsln $,  there exists a  $ \kqqm $--integer form (as
Hopf algebra)  $ {\Cal U}^{\scriptscriptstyle P} \big( {{\frak sl}(n+1)} \big)
$,  resp.~$ {\Cal U}^{\scriptscriptstyle P} \big( {\slgot(n+1)} \big) $
(cf.~[Ga], \S 3.4);  then we define  (cf.~[Ga], \S 4.3)
  $$  \eqalign {
   \calfQsln  &  := \left\{\, f \in \fqQsln \,\Big\vert\, \Big\langle f, {\Cal
U}^{\scriptscriptstyle P} \big( {\slgot (n+1)} \big) \Big\rangle \subseteq
\kqqm \,\right\}  \cr
   \calfPsln  &  := \left\{\, f \in \fqPsln \,\Big\vert\, \Big\langle f, {\Cal
U}^{\scriptscriptstyle Q} \big( {\slgot (n+1)} \big) \Big\rangle \subseteq
\kqqm \,\right\}  \cr }  $$
   \indent   From the very definition of  $ {\Cal U}^{\scriptscriptstyle P} $,
one sees that  $ \, \fqQtildesln \subseteq \calfQsln \, $,  and similarly,  $ \,
\fqPtildesln \subseteq \calfPsln \, $.  But even more, we shall prove in this
section that  $ \fqQtildesln $,  resp.~$ \fqPtildesln $,  is a "good enough
approximation" of  $ \calfQsln $,  resp.~$ \calfPsln $,  in the sense that they
have the same specialization at  $ q=1 $.
                                                \par
   One of the main points in [Ga] is the construction of a (topological) Hopf
algebra  $ U_q^{\scriptscriptstyle P} (\gerh) $,  with a  $ \kqqm $--integer
form $ \gerU^{\scriptscriptstyle P} (\gerh) $  which specializes to  $ \uh $
for  $ q \rightarrow 1 $.  The link with quantum function algebras is the
existence of an embedding of (topological) Hopf algebras
  $$  \xi_{\scriptscriptstyle P} \colon \, \fqPsln \llonghookrightarrow
U_q^{\scriptscriptstyle P} (\gerh) \; ;  $$
via this embedding one has  $ \, \xi_{\scriptscriptstyle P} \left( \calfPsln
\right) \subseteq \gerU^{\scriptscriptstyle P} (\gerh) \, $.  In addition, one
has also  $ \, \xi_{\scriptscriptstyle P} \left( \calfPsln
\right){\Big\vert}_{q=1} = \gerU^{\scriptscriptstyle P} (\gerh){\Big\vert}_{q=1}
\, $, so that  $ \, \gerU^{\scriptscriptstyle P} (\gerh) {\buildrel {q
\rightarrow 1} \over \llongrightarrow} \uh \, $  implies
  $$  \calfPsln {\buildrel {q \rightarrow 1} \over
\llongrightarrow} \uh \, .  $$
   \indent   The embedding  $ \, \xi_{\scriptscriptstyle P} \colon \, \fqPsln
\hookrightarrow U_q^{\scriptscriptstyle P} (\gerh) \, $  is the composition of
an embedding
  $$  \mu_{\scriptscriptstyle P} \colon \fqPg @>{\Delta}>> \fqPg \otimes \fqPg
@>{\rho_+ \otimes \rho_-}>> F_q^{\scriptscriptstyle P}[B_+] \otimes
F_q^{\scriptscriptstyle P}[B_-] @>{\vartheta_+ \otimes \vartheta_-}>>
{U_q^{\scriptscriptstyle P} (\gerb_-)}_{op} \otimes {U_q^{\scriptscriptstyle P}
(\gerb_+)}_{op}  $$
(where  $ \, G = SL(n+1) \, $,  and  $ B_\pm $  and  $ \gerb_\pm $  denotes
as usual Borel subgroups and Borel subalgebras) and an isomorphism
$ \nu_{\scriptscriptstyle P}^{\, -1} $  of a suitable subalgebra of
$ \, {U_q^{\scriptscriptstyle P}(\gerb_-)}_{op} \otimes
{U_q^{\scriptscriptstyle P} (\gerb_+)}_{op} \, $  (containing
$ \mu_{\scriptscriptstyle P} \! \left( \fqPsln \right) $)  with
$ U_q^{\scriptscriptstyle P} (\gerh) $;  hereafter,  $ H_{op} $  will denote
the (unique) Hopf algebra with the same structure of  $ H $  but for
comultiplication, which is turned into the opposite one.  Everything holds as
well with  $ P $  and  $ Q $  exchanging their roles; on
the other hand, since by definition is  $ \, \fqQsln \subseteq \fqPsln \, $,
$ \, U_q^{\scriptscriptstyle Q} (\gerh) \subseteq U_q^{\scriptscriptstyle P}
(\gerh) \, $,  and  $ \, \xi_{\scriptscriptstyle Q} = {\xi_{\scriptscriptstyle
P}}{\big\vert}_{\fqQsln} \, $,  $ \, \mu_{\scriptscriptstyle Q} =
{\mu_{\scriptscriptstyle P}}{\big\vert}_{\fqQsln} \, $,  it will be enough to
study  $ \mu_{\scriptscriptstyle P} \, $.  To this end, we have to revisit the
definition of  $ U_q^{\scriptscriptstyle M} (\slgot(n+1)) $  and its quantum
Borel subalgebras  $ U_q^{\scriptscriptstyle M}(\gerb_\pm) $  ($ \, M = Q, P
\, $), and the construction of quantum root vectors: this will be done in next
sections.  Here we recall the definition of
$ F_q^{\scriptscriptstyle P}[B_+] $  and  $ F_q^{\scriptscriptstyle P}[B_-] $
and the canonical epimorphisms  $ \rho_+ $  and  $ \rho_- \, $.
                                              \par
   $ F_q^{\scriptscriptstyle P}[B_+] $,
resp.~$ F_q^{\scriptscriptstyle P}[B_-] $,  is the unital associative
\hbox{$\kq$--alge}bra  generated by \,  $ \{ \rho_{ij} \mid i, j = 1, \ldots,
n+1; i \leq j \} $,  resp.~by \,  $ \{ \rho_{ij} \mid i, j = 1, \ldots, n+1; i
\geq j \} $,  \, with relations
  $$  \eqalignno {
   \rho_{ij} \rho_{ik} = q \, \rho_{ik} \rho_{ij} \; ,  \quad \quad  \rho_{ik}
\rho_{hk}  &  = q \, \rho_{hk} \rho_{ik}  &   \forall\, j<k, i<h  \qquad  \cr
   \rho_{il} \rho_{jk} = \rho_{jk} \rho_{il} \; ,  \quad \quad  \rho_{ik}
\rho_{jl} - \rho_{jl} \rho_{ik}  &  = \left( q - \qm \right) \, \rho_{il}
\rho_{jk}  \hskip90pt  &   \forall\, i<j, k<l  \qquad  \cr
   \rho_{1,1} \rho_{2,2} \cdots \rho_{n+1,n+1} = 1  &  {}  &  {}  \cr }  $$
for either algebras.  These are Hopf algebras too, with comultiplication given
by
  $$  \eqalign{
   \Delta (\rho_{ij})  &  = \sum_{k=i}^n \rho_{ik} \otimes \rho_{kj}  \qquad
\quad  \forall\, i, j  \qquad  \hbox{for
\ } F_q^{\scriptscriptstyle P}[B_+] \, ,  \cr
   \Delta (\rho_{ij})  &  = \sum_{k=1}^j \rho_{ik} \otimes \rho_{kj}  \qquad
\quad  \forall\, i, j  \qquad  \hbox{for
\ } F_q^{\scriptscriptstyle P}[B_-] \, ,  \cr } $$
counit by
  $$  \epsilon(\rho_{ij}) = \delta_{ij}  \qquad \quad  \forall\, i, j  $$
for either algebras, and antipode by
  $$  \eqalign{
   S(\rho_{ij})  &  = {(-q)}^{j-i} {det}_q^+ \left( {(\rho_{hk})}_{h \neq j}^{k
\neq i} \right)  \qquad \quad  \forall\, i, j  \qquad  \hbox{for \ }
F_q^{\scriptscriptstyle P}[B_+] \, ,  \cr
   S(\rho_{ij})  &  = {(-q)}^{j-i} {det}_q^- \left( {(\rho_{hk})}_{h \neq j}^{k
\neq i} \right)  \qquad \quad  \forall\, i, j  \qquad  \hbox{for \ }
F_q^{\scriptscriptstyle P}[B_-] \, ,  \cr }  $$
where  $ {det}_q^+ $,  resp.~$ {det}_q^- $,  is the expression that one gets
simply by setting  $ \, \rho_{ij} \equiv 0 \, $  for all  $ \, i>j \, $,
resp.~$ \, i<j \, $,  in  $ {det}_q \, $.
                                                  \par
   By the very definitions, two Hopf algebra epimorphisms exist
  $$  \rho_+ : \fqPsln \llongtwoheadrightarrow F_q^{\scriptscriptstyle P}[B_+]
\, ,  \qquad \quad  \rho_- : \fqPsln \llongtwoheadrightarrow
F_q^{\scriptscriptstyle P}[B_-] \, ,  $$
given by
  $$  \eqalign{
   \rho_+ : \quad \rho_{ij} \mapsto \rho_{ij} \; \quad \forall\, i \leq j \, ,
&  {}  \; \qquad \rho_{ij} \mapsto 0 \; \; \forall\, i > j  \cr
   \rho_- : \quad \rho_{ij} \mapsto \rho_{ij} \; \quad \forall\, i \geq j \, ,
&  {}  \; \qquad \rho_{ij} \mapsto 0 \; \; \forall\, i < j \, .  \cr }  $$

\vskip7pt

  {\bf 5.2  The quantum algebras  $ \uqgln $,  $ \uqsln $,  and
$ U_q^{\scriptscriptstyle P} \left( \gerb_\pm \right) $.} \  We
recall (cf.~for instance [GL]) the definition of the quantized universal
enveloping algebra  $ \uqgln $:  it is the associative algebra with 1 over
$ \kq $  with generators
  $$  F_1, F_2, \dots, F_n, G_1^{\pm 1}, G_2^{\pm 1}, \dots, G_n^{\pm 1},
G_{n+1}^{\pm 1}, E_1, E_2, \dots, E_n  $$
and relations
  $$  \displaylines {
   \hfill   G_i G_i^{-1} = 1 = G_i^{-1} G_i \, ,  \qquad  G_i^{\pm 1} G_j^{\pm
1} = G_j^{\pm 1} G_i^{\pm 1}  \qquad   \hfill  \forall \, i, j \qquad  \cr
   \hfill   G_i F_j G_i^{-1} = q^{\delta_{i,j+1} - \delta_{i,j}} F_j \, ,
\qquad  G_i E_j G_i^{-1} = q^{\delta_{i,j} - \delta_{i,j+1}} E_j  \qquad
\hfill  \forall \, i, j \qquad  \cr
   \hfill   E_i F_j - F_j E_i = \delta_{i,j} {{\; G_i G_{i+1}^{-1} - G_i^{-1}
G_{i+1} \;} \over {\; q - \qm \;}}  \qquad \qquad   \hfill  \forall \, i, j
\qquad  \cr
   \hfill   \quad  E_i E_j = E_j E_i \, ,  \qquad  F_i F_j = F_j F_i   \hfill
\forall \; i, j \colon \vert i - j \vert > 1 \, \phantom{.} \;  \cr
   \hfill   E_i^2 E_j - \left( q + \qm \right) E_i E_j E_i + E_j E_i^2 = 0
\hfill  \forall \; i, j \colon \vert i - j \vert = 1 \, \phantom{.} \;  \cr
   \hfill   F_i^2 F_j - \left( q + \qm \right) F_i F_j F_i + F_j F_i^2 = 0
\hfill  \forall \; i, j \colon \vert i - j \vert = 1 \, . \;  \cr }  $$
   \indent   Moreover,  $ \uqgln $  has a Hopf algebra structure, given by
  $$  \displaylines {
   \hfill   \Delta \left( F_i \right) = F_i \otimes G_i^{-1} G_{i+1} + 1 \otimes
F_i \, ,  \; \qquad  S \left( F_i \right) = - F_i G_i G_{i+1}^{-1}
\, , \; \qquad  \epsilon \left( F_i \right) = 0   \hfill  \forall \, i \,
\phantom{.} \;  \cr
   \hfill   \qquad  \Delta \left( G_i^{\pm 1} \right) = G_i^{\pm 1} \otimes
G_i^{\pm 1} \, ,  \qquad \qquad \qquad  S \left( G_i^{\pm 1} \right) = G_i^{\mp
1} \, ,  \quad \qquad  \epsilon \left( G_i^{\pm 1} \right) = 1   \hfill  \forall
\, i \, \phantom{.} \;  \cr
   \hfill   \Delta \left( E_i \right) = E_i \otimes 1 + G_i G_{i+1}^{-1} \otimes
E_i \, ,  \qquad  S \left( E_i \right) = - G_i^{-1} G_{i+1} E_i \, ,  \qquad
\epsilon \left( E_i \right) = 0   \hfill  \forall \, i \, . \;
\cr }  $$
   The algebras  $ \uqPsln $  and  $ \uqQsln $   --- defined as in  [Ga], \S 3
---   can be realized as Hopf subalgebras or quotients of  $ \uqgln $.  Namely,
define elements
  $$  L_i := G_1 \cdots G_i \, ,  \; \quad L_i^{-1} := G_1^{-1} \cdots G_i^{-1}
\, ,  \; \quad K_i := G_i G_{i+1}^{-1} \, ,  \; \quad K_i^{-1} := G_i^{-1}
G_{i+1}  $$
for all  $ \, i=1,\dots,n \, $.  Then  $ L_{n+1} $  is a central element of
$ \uqgln $,  and  $ \uqPsln $  is (isomorphic to) the subalgebra of  $ \uqgln $
generated by  $ \{ F_1, \dots, F_n, L_1^{\pm 1}, \dots,
                             L_n^{\pm 1}, $\break
$ E_1, \dots, E_n \} $   --- this corresponds to  $ \, \slgot(n+1)
\hookrightarrow \gl(n+1) \, $  ---   and to the quotient of  $ \uqgln $  modulo
the ideal (which is a  {\sl Hopf} ideal) generated by  $ \left( L_{n+1} - 1
\right) $   --- which corresponds to  $ \, \slgot(n+1) \cong \gl(n+1) \big/ \!
\left( I_{n+1} \right) = \gl(n+1) \big/ Z\big(\gl(n+1)\big) \, $.  Similarly,
the algebra  $ \uqQsln $  is (isomorphic to) the subalgebra of  $ \uqgln $
generated by  $ \{ F_1, \dots, F_n, K_1^{\pm 1}, \dots, K_n^{\pm 1}, E_1, \dots,
E_n \} $   --- this again corresponds to  $ \, \slgot(n+1) \hookrightarrow
\gl(n+1) \, $.
                                                   \par
   A last word about "toral elements"  $ G_i $,  $ L_i $,  $ K_i $.  In the
"classical" framework we have toral elements  $ \, h_i = h_{\alpha_i} \, $  in
the (diagonal) Cartan subalgebra of  $ \, \slgot(n+1) \big( \subseteq \gl(n+1)
\big) \, $,  given by  $ \, h_i = M_{i,i} - M_{i+1,i+1} \, $  (where
$ M_{r,s} $  denotes a square matrix of size  $ n+1 $  as in \S 2); similarly,
letting  $ \, \ell_i = M_{1{}1} + M_{2{}2} + \cdots + M_{i{}i} \, $,  we have
$ \, h_i = - \ell_{i-1} + 2 \ell_{i} - \ell_{i+1} \, $.  Now,  $ G_i $  is the
$ q $--analogue  of  $ M_{i,i} $   --- in fact, on the standard representation
it acts exactly as  $ \exp \big ( h M_{i,i} \big) $, where  $ \, h := \log (q)
\, $  ---   therefore  $ \, K_i := G_i G_{i+1}^{-1} = \exp \! \big ( h \,
M_{i,i} \big) \cdot {\exp \! \big ( h \, M_{i+1,i+1} \big)}^{-1} = \exp \! \big
( h \, (M_{i,i} - M_{i+1,i+1}) \big) = \exp \! \big ( h \cdot h_i \big) \, $  is
exactly the  $ q $--analogue  of $ h_i $;  similarly,  $ L_i $  is the
$ q $--analogue  of  $ \ell_i $.
                                                   \par
   As for quantum Borel subalgebras, we recall that
$ U_q^{\scriptscriptstyle P} (\gerb_+) $,
resp.~$ U_q^{\scriptscriptstyle P} (\gerb_-) $,  is   --- by definition ---
the subalgebra of  $ \uqPsln $  generated by  $ \, \{ L_1, \dots, L_n \} \cup \{
E_1, \dots, E_n \} \, $,  resp.~by  $ \, \{ L_1, \dots, L_n \} \cup \{ F_1,
\dots, F_n \} \, $;  similar definitions occur for  $ U_q^{\scriptscriptstyle Q}
(\gerb_\pm) $,  with  $ K_i $'s  instead of  $ L_i $'s.  All these are in fact
{\sl Hopf}  subalgebras.
                                                  \par
   Finally, there exist Hopf algebra isomorphisms
  $$  \vartheta_+ : F_q^{\scriptscriptstyle P} [B_+] {\buildrel \cong \over
\llongrightarrow} {U_q^{\scriptscriptstyle P} (\gerb_-)}_{op} \, ,  \qquad
\vartheta_- : F_q^{\scriptscriptstyle P} [B_-] {\buildrel \cong \over
\llongrightarrow} {U_q^{\scriptscriptstyle P} (\gerb_+)}_{op}  $$
which are uniquely determined by
  $$  \displaylines{
   \hfill   \vartheta_+ \left( \rho_{ii} \right) := L_{i-1} L_i^{-1} = G_i^{-1}
\; ,  \quad  \vartheta_+ \left( \rho_{i,i+1} \right) := - \fbar_i L_i
L_{i+1}^{-1} = - \fbar_i G_{i+1}^{-1} \; ,   \hfill  \forall\; i \;  \cr
   \hfill   \vartheta_- \left( \rho_{ii} \right) := L_{i-1}^{-1} L_i = G_i
\; ,  \quad  \vartheta_- \left( \rho_{i+1,i} \right) := + L_i^{-1} L_{i+1}
\ebar_i = + G_{i+1} \ebar_i \; ,   \hfill  \forall\; i \;  \cr } $$
(here we set  $ \, L_0 := 1 \, $,  $ \, L_{n+1} := 1 \, $),  where  $ \, \fbar_i
:= \left( q - \qm \right) F_i \, $,  $ \, \ebar_i := \left( q - \qm \right) E_i
\, $.

\vskip7pt

  {\bf 5.3  Quantum root vectors.} \  Quantum root vectors are essential in
[Ga]: according to a general recipe provided by Lusztig, they are constructed by
means of braid group operators  $ T_i \, (i=1, \dots, n) $,  which in our case
are given by (with the normalization of [Ga])
  $$  T_i \colon
   \cases
      F_j \mapsto -K_j^{-1} E_j \, ,  \hskip27pt  \quad  K_j \mapsto K_j^{-1}
\, ,  \hskip5pt \quad  E_j \mapsto -F_j K_j \, ,  &  \;\;  \hbox{if}  \;\;\;
\vert i-j \vert = 0  \\
      F_j \mapsto F_j \, ,   \hskip55,4pt  \quad  K_j \mapsto K_j \, ,
\hskip12,1pt  \quad  E_j \mapsto E_j \, ,  &  \;\;  \hbox{if}  \;\;\; \vert i-j
\vert > 1  \\
      F_j \mapsto -F_j F_i + q \hskip1pt F_i F_j \, ,  \quad \hskip-1pt  K_j
\mapsto K_i K_j \, ,  \quad  E_j \mapsto -E_i E_j + \qm E_j E_i \, ,  &  \;\;
\hbox{if}  \;\;\; \vert i-j \vert = 1  \\
   \endcases  $$
   \indent   Thus letting  $ \; {[x,y]}_q := x y - q \, y x \; $  be the
$ q $--bracket  of  $ x $  and  $ y $  we have
  $$  T_i (F_j) = q \, {[F_i,F_j]}_{\qm} = -{[F_j,F_i]}_q \; ,  \; \;  T_i (E_j)
= -{[E_i,E_j]}_{\qm} = \qm {[E_j,E_i]}_q  \; \;  \hbox{\, for \,} \; \vert i-j
\vert = 1 .  $$
   \indent   Consider now the case  $ \, G = SL(n+1) \, $,  $ \, \gerg =
\slgot(n+1) \, $.  We want to compare Lusztig's construction of quantum root
vectors with another one (which is used, for instance, in [Ji] and in [Ta]).
                                                   \par
   In the Lie algebra  $ \slgot(n+1) $  we have matrices  $ M_{ij} $  ($ i, j
\in \{1,2,\dots,n\}, i<j $)  such that
  $$  M_{i,i+1} = e_i \; ,  \qquad  [M_{ik}, M_{kj}] = M_{ij} \; ,  \qquad
\quad  \forall \, i<k<j  $$
where  $ e_i $  denotes the  $ i $--th  Chevalley generator of the positive part
of  $ \slgot(n+1) $,  the Lie subalgebra  $ \gern_+ $  of strictly upper
triangular matrices.  Moreover, these matrices are root vectors, i.~e.~weight
vectors   --- for the adjoint action ---   with roots as weights: namely,  $
M_{ij} $  has weight the positive root  $ \, \alpha(i,j) := \varepsilon_i -
\varepsilon_j \, $,  where  $ \varepsilon_k $  denotes the  $ (k,k) $--th
coordinate function on matrices.  {\it Notice that the same is true for}
$ \, - [M_{ik}, M_{kj}] = -M_{ij} \, $.  Then the  {\sl level}  of  $ M_{ij} $,
defined to be  $ j - i $,  is the  {\sl height}  of the root  $ \alpha(i,j) $.
A similar situation occurs for the subalgebra  $ \gern_- $  of strictly lower
triangular matrices, with elements  $ M_{ij} $  ($ i>j $)  as root vectors of
weight the negative roots  $ \alpha(i,j) $  of height  $ j-i $.  In particular,
$ \big\{\, M_{ij} \,\big\vert\, i<j \,\big\} $  is a basis of  $ \gern_+ $,
and  $ \big\{\, M_{ij} \,\big\vert\, i>j \,\big\} $  is a basis of  $ \gern_- $.
                                                 \par
  The situation described above can be quantized.  In fact in the standard
representation of  $ \uqPsln $  or  $ \uqQsln $  we still have  $ \, E_i \mapsto
M_{i,i+1} \, $,  i.e.~$ E_i $  acts as the matrix  $ M_{i,i+1} $;  then we
define (for all  $ \, i < j-1 \, $)
  $$  \eqalign{
   E_{i,i+1} := E_i \; ,  &  \qquad  E_{i,j} := -{\left[ E_{i,j-1}, E_{j-1,j}
\right]}_{\qm} = \qm {\left[ E_{j-1,j}, E_{i,j-1} \right]}_q  \cr
   F_{j+1,j} := F_i \; ,  &  \qquad  F_{j,i} := q \, {\left[ F_{j-1,i},
F_{j,j-1} \right]}_{\qm} = - {\left[ F_{j,j-1}, F_{j-1,i} \right]}_q  \cr }
\eqno (5.1)  $$
(notice the occurrence of the  "$ - $"  sign, which does not appear in [Ji], nor
in [Ta]); then remark also that in the standard representation of  $ \uqPsln $
or  $ \uqQsln $  we have  $ \, E_{ij} \mapsto {(-1)}^{j-i-1} M_{ij} \, $  and
$ \, F_{ji} \mapsto {(-1)}^{j-i-1} M_{ji} \, $,  for all  $ \, i < j \, $.
                                                   \par
   Now look at roots, for instance positive ones: they form the set
  $$  R^+ := \big\{\, \alpha(i,j) \,\big\vert\, i, j= 1, \dots, n+1; i<j
\,\big\} \; ;  $$
if we set  $ \, n(i,j) := j - i + \sum_{h=1}^{i-1} (n-h) \, $,  we obtain a
total ordering of  $ R^+ $  by
   $$  \alpha(i,j) \preceq \alpha(h,k) \; \Longleftrightarrow \; n(i,j) \leq
n(h,k)  $$
so that  $ \, R^+ = \left\{ \alpha^1, \alpha^2, \dots, \alpha^N \right\} \, $,
with  $ \, \alpha(i,j) =: \alpha^{n(i,j)} \, $,  $ \, N:= {n+1 \choose 2} \, $.
                                                   \par
   The first key point is the following lemma, whose proof is trivial:

\vskip5pt

\proclaim{Lemma}  The previously defined ordering of
$ R^+ $  is convex, that is
  $$  \alpha, \beta, \alpha + \beta \in R^+, \alpha \preceq \beta \;
\Longrightarrow \; \alpha \preceq \alpha +\beta \preceq \beta \, .   \, \;
\square  $$
\endproclaim

\vskip5pt

   By a theorem of Papi (cf.~[Pa]) we know that every total ordering of
$ R^+ $  which is convex is associated to a unique minimal expression of  $ \,
w_0 = (n \, n-1 \, \dots \, 3 \,2 \,1) \, $,  the longest element of the Weyl
group  $ S_n $  of  $ \slgot(n+1) $:  this means that,  if  $ \, w_0 = s_{i_1}
s_{i_2} \cdots s_{i_{N-1}} s_{i_N} \, $  is a minimal expression of  $ w_0 $,
then the ordering of  $ R^+ $  is given by
  $$  \alpha^1 = \alpha_{i_1} \, ,  \, \alpha^2 = s_{i_1} \big( \alpha_{i_2}
\big) \, ,  \, \dots \, ,  \, \alpha^k = s_{i_1} s_{i_2} \cdots s_{i_{k-1}}
\big( \alpha_{i_k} \big) \, ,  \, \dots \, ,  \, \alpha^N = s_{i_1} s_{i_2}
\cdots s_{i_{N-1}} \big( \alpha_{i_N} \big) \, .  $$
                                                 \par
   In the present situation we can also write down explicitely the minimal
expression of  $ w_0 $  afforded by the given ordering of   $ R^+ \, $:  it is
  $$  w_0 = s_1 s_2 s_3 \cdots s_{n-1} s_n s_1 s_2 \cdots s_{n-1} s_1 s_2 \cdots
s_{n-3} s_1 s_2 \cdots s_4 s_1 s_2 s_3 s_1 s_2 s_1 \; .  $$
   \indent   Starting from any minimal expression of
$ w_0 $,  Lusztig constructs quantum root vectors via the formulae
  $$  E_{\alpha^k} := T_{w_{k-1}} \! \left( E_{i_k} \right) = T_{i_1} T_{i_2}
\cdots T_{i_{k-1}} \! \left( E_{i_k} \right) \, ,  \qquad  F_{\alpha^k} :=
T_{w_{k-1}} \! \left( F_{i_k} \right) = T_{i_1} T_{i_2} \cdots T_{i_{k-1}} \!
\left( F_{i_k} \right)  $$
where  $ \, w_0 = s_{i_1} s_{i_2} \cdots s_{i_N} \, $  is the given minimal
expression,  $ \, w_h := s_{i_1} s_{i_2} \cdots s_{i_h} \, $,  and  $ \,
\alpha^k = w_{k-1} \left( \alpha_{i_k} \right) \, $  gives the associated convex
ordering of  $ R^+ $.  In our case this construction gives quantum root vectors
(for all  $ \, i < j \, $)
  $$  \displaylines{
   E_{\alpha^{n(i,j)}} := T_1 T_2 \cdots T_n T_1 T_2 \cdots T_{n-1} T_1
\cdots T_{n-i+1} T_1 T_2 \cdots T_{j-i-1} \left( E_{j-i} \right) \, ,  \cr
   F_{\alpha^{n(i,j)}} := T_1 T_2 \cdots T_n T_1 T_2 \cdots T_{n-1} T_1
\cdots T_{n-i+1} T_1 T_2 \cdots T_{j-i-1} \left( F_{j-i} \right) \, .  \cr }  $$

\vskip5pt

\proclaim{Theorem}  For all  $ \, i, j= 1, 2, \dots, n+1 \, $  with  $ \, i < j
\, $  we have
  $$  E_{\alpha^{n(i,j)}} = E_{ij} \; ,  \quad \qquad  F_{\alpha^{n(i,j)}} =
F_{ji} \; .  $$
\endproclaim

\demo{Proof}  We make the proof for the  "E" case, the "F" case is the like.
                                                 \par
   We have two possibilities,  $ \, j=i+1 \, $  and  $ \, j>i+1 \, $.  If  $ \,
j=i+1 \, $  we have
  $$  \alpha^{n(i,j)} = s_1 s_2 \cdots s_n s_1 s_2 \cdots s_{n-i+2} s_1 s_2
\cdots s_{n-i+1} (\alpha_1) = s_1 s_2 \cdots s_n s_1 s_2 \cdots s_{n-i+2} s_1
s_2 (\alpha_1)  $$
because  $ \, s_h (\alpha_k) = \alpha_k \, $  for all  $ \, \vert h-k \vert > 1
\, $;  since  $ \, s_1 s_2 (\alpha_1) = \alpha_2 \, $,  we get
  $$  \displaylines{
   {\ } \quad   \alpha^{n(i,j)} = s_1 s_2 \cdots s_n s_1 s_2 \cdots s_{n-i+2}
(\alpha_2) =   \hfill {\ }  \cr
   {\ } \hfill   = s_1 s_2 \cdots s_n s_1 s_2 \cdots s_{n-i+3} s_1 s_2 s_3
(\alpha_2) = s_1 s_2 \cdots s_n s_1 s_2 \cdots s_{n-i+3} (\alpha_3) \; ;
\cr }  $$
thus iterating we obtain
  $$  \alpha^{n(i,j)} = \alpha_i  \; ;  $$
then by Lusztig's work we know that  $ \, E_{\alpha^{n(i,j)}} = E_{\alpha_i} =
E_i \, $;  on the other hand,  $ \, j=i+1 \, $  gives  $ \, E_{i,j} = E_{i,i+1}
= E_i \, $  by definition, hence  $ \, E_{\alpha^{n(i,j)}} = E_{i,j} \, $,
q.e.d.
                                                 \par
   If  $ \, j>i+1 \, $  we apply the definitions to get
  $$  \displaylines{
   {\ }  E_{\alpha^{n(i,j)}} := T_1 T_2 \cdots T_n T_1 T_2 \cdots T_{n-1} T_1
\cdots T_{n-i+1} T_1 T_2 \cdots T_{j-i-2} T_{j-i-1} \left( E_{j-i} \right)
\hfill  \cr
   = T_1 T_2 \cdots T_n T_1 T_2 \cdots T_{n-1} T_1 \cdots T_{n-i+1} T_1 T_2
\cdots T_{j-i-2} \left( - {\left[ E_{j-i-1}, E_{j-i} \right]}_{\qm} \right) =
\cr
   = - {\left[ T_1 \cdots T_{n-i+1} T_1 \cdots T_{j-i-2} \left( E_{j-i-1}
\right), T_1 \cdots T_{n-i+1} T_1 \cdots T_{j-i-2} \left( E_{j-i} \right)
\right]}_{\qm} =  \cr
   = - {\left[ E_{\alpha^{n(i,j-1)}} \, , T_1 T_2 \cdots T_n T_1 T_2 \cdots
T_{n-i+1} T_1 T_2 \cdots T_{j-i-2} \left( E_{j-i} \right) \right]}_{\qm} \; ;
\cr }  $$
but now we have
  $$  \displaylines{
   {\ } T_1 T_2 \cdots T_n T_1 T_2 \cdots T_{n-i+2} \cdots T_1 T_2 \cdots
T_{n-i+1} T_1 T_2 \cdots T_{j-i-2} \left( E_{j-i} \right) =   \hfill {\ }  \cr
   = T_1 T_2 \cdots T_n T_1 T_2 \cdots T_{n-i+2} T_1 T_2 \cdots T_{j-i-1}
T_{j-i} T_{j-i+1} \left( E_{j-i} \right) =  \cr
   {\ } \hfill   = T_1 T_2 \cdots T_n T_1 T_2 \cdots T_{n-i+2} T_1 T_2 \cdots
T_{j-i-1} \left( E_{j-i+1} \right)  \cr }  $$
for  $ \, T_k T_{k+1} \left( E_k \right) = E_{k+1} \, $  for all  $ k \, $;
then a simple iteration gives
  $$  T_1 T_2 \cdots T_n T_1 T_2 \cdots T_{n-i+2} \cdots T_1 T_2 \cdots
T_{n-i+1} T_1 T_2 \cdots T_{j-i-2} \left( E_{j-i} \right) = E_{j-1} \; .  $$
   \indent   Therefore we have
  $$  E_{\alpha^{n(i,j)}} = - {\left[ E_{\alpha^{n(i,j-1)}} \, , E_{j-1}
\right]}_{\qm} \; ;  $$
since  $ \, (j-1) - i < j - i \, $,  we can use induction on the height of roots
and assume  $ \, E_{\alpha^{n(i,j-1)}} = E_{i,j-1} \, $;  on the other hand,  $
\, E_{j-1} = E_{j-1,j} \, $,  by definition; thus we find  $ \,
E_{\alpha^{n(i,j)}} = - {\left[ E_{i,j-1} \, , E_{j-1,j} \right]}_{\qm} \, $,
and   --- by definition again ---   the right-hand-side term is nothing but  $
E_{i,j} \, $.  The claim follows.   $ \square $
\enddemo

\vskip7pt

  {\bf 5.4  The embedding  $ \, \mu_{\scriptscriptstyle M} \colon \fqPsln
\llonghookrightarrow {U_q^{\scriptscriptstyle P} (\gerb_-)}_{op} \otimes
{U_q^{\scriptscriptstyle P} (\gerb_+)}_{op} \, $  and specialization results.}
\  We are now ready to go on with the analysis of the embedding
  $$  \mu_{\scriptscriptstyle P} \colon \fqPg @>{\Delta}>> \fqPg \otimes \fqPg
@>{\rho_+ \otimes \rho_-}>> F_q^{\scriptscriptstyle P}[B_+] \otimes
F_q^{\scriptscriptstyle P}[B_-] @>{\vartheta_+ \otimes \vartheta_-}>>
{U_q^{\scriptscriptstyle P} (\gerb_-)}_{op} \otimes {U_q^{\scriptscriptstyle P}
(\gerb_+)}_{op}  $$
($ \, G = SL(n+1) \, $):  here  $ \Delta $  is the comultiplication of the Hopf
algebra  $ F_q^{\scriptscriptstyle M}[SL(n+1)] $,  $ \, \rho_\pm \colon
F_q^{\scriptscriptstyle M} [SL(n+1)] \llongtwoheadrightarrow
F_q^{\scriptscriptstyle M}[B_\pm] \, $  are the natural epimorphisms of \S 5.1,
and the maps  $ \, \vartheta_\pm \colon F_q^{\scriptscriptstyle M}[B_\pm]
{\buildrel \cong \over \llongrightarrow} {U_q^{\scriptscriptstyle M}
(\gerb_\mp)}_{op} \, $  are the isomorphisms given in \S 5.2.
                                                \par
  To begin with, we go back to the study of  $ \fqQtildesln $:  lifting the
identities (4.2) gives relations (inside  $ \fqQtildesln \, $)
  $$  \hbox{ $ \eqalign{
   r_{ij}  &  = + \big[ r_{i,j-1} \, , \, r_{j-1,j} \big] + {\Cal O}(q-1) = -
\big[ r_{j-1,j} \, , r_{i,j-1} \big] + {\Cal O}(q-1)  \quad \;  \forall\, i <
j-1 \, ,  \cr
   r_{ji}  &  = - \big[ r_{j,j-1} \, , \, r_{j-1,i} \big] + {\Cal O}(q-1) = +
\big[ r_{j-1,i} \, , \, r_{j,j-1} \big] + {\Cal O}(q-1)  \quad \;  \forall\, j >
i+1 \, .  \cr } $ }   \eqno (5.1)  $$
   \indent   Now look at the isomorphism  $ \vartheta_+ \colon
F_q^{\scriptscriptstyle P}[B_+] \, {\buildrel \cong \over \llongrightarrow} \,
{U_q^{\scriptscriptstyle P} (\gerb_-)}_{op} $;  it maps  $ \, \rho_{ii} = r_{ii}
\, $
                    onto\break
$ \, L_{i-1} L_i^{-1} = G_i^{-1} \, $,  and  $ \rho_{i,i+1} $  onto  $ \, -
\fbar_i L_i L_{i+1}^{-1} = - \fbar_i + G_{i+1}^{\,-1} \, $,  hence
$ r_{i,i+1} $
          onto
          %\break
$ \, - F_i L_i L_{i+1}^{-1} \, $.  But remark that  $ \, r_{ii} = 1 + {\Cal
O}(q-1) \, $,  and so also  $ \, G_i = 1 + {\Cal O}(q-1) \, $,  $ \, K_i = 1 +
{\Cal O}(q-1) \, $,  $ \, L_i = 1 + {\Cal O}(q-1) \, $,  where the symbol
$ {\Cal O}(q-1) $  denotes some element of  $ \, (q-1) \,
\gerU^{\scriptscriptstyle P} (\gerb_-) \, $  (notations of [Ga]):  thus we have
also  $ \, \vartheta_+ \big( r_{i,i+1} \big) = - F_i + {\Cal O}(q-1) \, $.
                                        \par
   Therefore, using the first relation of (5.1), Theorem 5.3, and the fact that
$ \, {[\ ,\ ]}_q \equiv {[\ ,\ ]}_{q^{-1}} \equiv [\ ,\ ] \mod\, (q-1) \, $,  by
a simple iterating procedure we find that
  $$  \vartheta_+ \big( r_{ij} \big) = {(-1)}^{j-i} F_{j,i} + {\Cal O}(q-1) =
{(-1)}^{j-i} F_{\alpha^{n(i,j)}} + {\Cal O}(q-1)  \qquad \quad \forall \; i < j
\, .  $$
   Similarly, by the same arguments we can prove that
  $$  \vartheta_- \big( r_{ji} \big) = {(-1)}^{j-i+1} E_{i,j} + {\Cal O}(q-1) =
{(-1)}^{j-i+1} E_{\alpha^{n(i,j)}} + {\Cal O}(q-1)  \qquad \quad  \forall \; j >
i \, .  $$
   \indent   Now, from the definition of  $ \Delta $  and  $ \vartheta_\pm $  we
get a description of  $ \, \mu_{\scriptscriptstyle P}( r_{i,i+1}) \, $  as
follows:
  $$  \displaylines{
   \mu_{\scriptscriptstyle P} (r_{i,i+1}) = (\vartheta_+ \otimes \vartheta_-)
\bigg( (\rho_+ \otimes \rho_-) \Big( \Delta \big( r_{i,i+1} \big)
\Big) \bigg) =   \hfill {\ }  \cr
   = (\vartheta_+ \otimes \vartheta_-) \Bigg( (\rho_+ \otimes \rho_-) \bigg(
r_{i,i} \otimes r_{i,i+1} + r_{i,i+1} \otimes r_{i+1,i+1} + \left( q - \qm
\right) \sum_{\Sb  k = 1 \\  k \neq i,j \\  \endSb}^{n+1} r_{i,k} \otimes
r_{k,i+1} \bigg) \Bigg) =  \cr
   = (\vartheta_+ \otimes \vartheta_-) \left( r_{i,i+1} \otimes r_{i+1,i+1} +
\left( q - \qm \right) \sum_{k=i+2}^{n+1} r_{i,k} \otimes r_{k,i+1} \right) =
\cr
   = - F_i \,G_{i+1}^{\,-1} \otimes G_{i+1} + {\Cal O}(q-1) + \left( q - \qm
\right) \sum_{k=i+2}^{n+1} F_{\alpha(k,i)} \otimes E_{\alpha(i+1,k)} + {\Cal
O}{\left( {(q-1)}^2 \right)} =  \cr
   {\ } \hfill   = - F_i \, G_{i+1}^{\,-1} \otimes G_{i+1} + {\Cal
O}(q-1) = - F_i \otimes 1 + {\Cal O}(q-1)  \cr }  $$
where the symbol  $ {\Cal O}(q-1) $  denotes some element of  $ \, (q-1) \,
\gerU^{\scriptscriptstyle P} (\gerb_-) \otimes \gerU^{\scriptscriptstyle P}
(\gerb_+) \, $  (note that the last algebra is mapped by
$ \nu_{\scriptscriptstyle P}^{\, -1} $  into  $ \gerU^{\scriptscriptstyle P}
(\gerh) \, $);  recalling from [Ga] the definition of  $ \,
\nu_{\scriptscriptstyle P}^{\, -1} : \gerU^{\scriptscriptstyle P} (\gerb_-)
\otimes \gerU^{\scriptscriptstyle P} (\gerb_+) \longrightarrow
\gerU^{\scriptscriptstyle P} (\gerh) \, $,  for  $ \, \xi_{\scriptscriptstyle P}
= \mu_{\scriptscriptstyle P} \circ \nu_{\scriptscriptstyle P}^{\, -1} \, $  we
find
  $$  \xi_{\scriptscriptstyle P} \left( r_{i,i+1} \right) = - F_i + {\Cal
O}(q-1)  $$
and in general  $ \, \xi_{\scriptscriptstyle P} \left( r_{i,j} \right) =
{(-1)}^{j-i} F_{j,i} + {\Cal O}(q-1) \, $,  for all  $ \, i < j \, $.  A similar
analysis yields
  $$  \xi_{\scriptscriptstyle P} \left( r_{i+1,i} \right) = + E_i + {\Cal
O}(q-1)  $$
and in general  $ \, \xi_{\scriptscriptstyle P} \left( r_{j,i} \right) =
{(-1)}^{j-i+1} E_{i,j} + {\Cal O}(q-1) \, $,  for all  $ \, j > i \, $,  and
also
  $$  \xi_{\scriptscriptstyle P} \left( r_{i,i} \right) = G_i + {\Cal O} \left(
{(q-1)}^2 \right) = L_{i-1}^{-1} L_i + {\Cal O} \left( {(q-1)}^2 \right)
\, ;  $$
furthermore, the latter implies
  $$  \displaylines{
   {} \quad   \xi_{\scriptscriptstyle P} \left( \varphi_i \right) =
\xi_{\scriptscriptstyle P} \left( { \; r_{i,i} - r_{i+1,i+1} \; \over \; q - 1
\; } \right)
= { \; G_i - G_{i+1} \; \over \; q - 1 \; } + {\Cal O}(q-1) =   \hfill {\ }  \cr
   = G_{i+1} \cdot { \; G_i G_{i+1}^{-1} - 1 \; \over \; q - 1 \; } + {\Cal
O}(q-1) =
G_{i+1} \cdot { \; K_i - 1 \; \over \; q - 1 \; } + {\Cal O}(q-1) =  \cr
   {\ } \hfill   = { \; K_i - 1 \; \over \;
q - 1 \; } + {\Cal O}(q-1) = \left({ K_i ; 0 \atop 1} \right) + {\Cal O}(q-1) \,
,  \cr
   {} \quad   \xi_{\scriptscriptstyle P} \left( \psi_i \right) =
\xi_{\scriptscriptstyle P} \left( { \; r_{1,1} \cdots r_{i,i} - 1 \; \over \; q
- 1 \; } \right) = { \; G_1 \cdots G_i - 1 \; \over \; q - 1 \; } + {\Cal
O}(q-1) =   \hfill {\ }  \cr
{\ } \hfill   = { \; L_i - 1 \; \over \; q - 1 \; } + {\Cal O}(q-1) = \left({
L_i ; 0 \atop 1} \right) + {\Cal O}(q-1) \, ,  \cr
   {} \quad   \xi_{\scriptscriptstyle P} \left( \chi_i \right) =
\xi_{\scriptscriptstyle P} \left( { \; r_{i,i} - 1 \; \over \; q - 1 \; }
\right) = { \; G_i - 1 \; \over \; q - 1 \; } + {\Cal O}(q-1) = \left({ G_i ; 0
\atop 1} \right) + {\Cal O}(q-1) \, .  \cr }  $$

\vskip5pt

   {\bf Remark:} \  in light of the previous formulae and of the remarks in \S
5.2 about toral elements of  $ \uqgln $, one has that
  $$  \displaylines{
   \varphi_i \, \hbox{\ is the  $ q $--analogue  of \ } h_i \, ,  \cr
   \psi_i \, \hbox{\ is the  $ q $--analogue  of \ } \ell_i \, ,  \cr
   \chi_i \, \hbox{\ is the  $ q $--analogue  of \ } M_{i,i} \, ;  \cr }  $$
on the other hand, the definition of these elements can also be motivated
directly in quantum matrix terms.  In the classical framework, we have
  $$  h_i = M_{i,i} - M_{i+1,i+1} \; ,  \qquad  \ell_i = M_{1,1} + M_{2,2} +
\cdots + M_{i,i} \; ;   \eqno (5.2)  $$
   \indent   Now,  $ r_{i,i} $  is the  $ i $--th  diagonal quantum matrix
coefficient, that is the  $ q $--analogue  of the "classical" matrix
coefficient  $ M_{i,i} $;  hence we should have, in principle,
  $$  {\; r_{i,i} - 1 \; \over \; q-1 \;}{\Big\vert}_{q=1} = M_{i,i} \; ;  $$
such a relation is completely meaningful, and was our reason to define  $ \,
\chi := {\; r_{i,i} - 1 \; \over \; q-1 \;} \, $.  As  $ G_i $  too is a
$ q $--analogue  of  $ M_{i,i} $,  the relation
  $$  \xi_{\scriptscriptstyle P} \left( \chi_i \right) = \left( G_i ; \, 0 \atop
1 \right) + {\Cal O}(q-1)  $$
is not surprising.  By the way, we remark that the special relation
  $$  \displaylines{
   {} \quad   \sum_{i=1}^{n+1} r_{1,1} r_{2,2} \cdots r_{i-1,i-1} \chi_i =
\hfill {\ }  \cr
   {\ } \hfill   = \sum_{\sigma \in S_{n+1} \setminus \{1\}} {(-q)^{l(\sigma)}}
{(q-1)}^{e(\sigma)-1} {\left( 1 + \qm \right)}^{e(\sigma)} r_{1,\sigma(1)}
r_{2,\sigma(2)} \cdots r_{n+1,\sigma(n+1)}  \cr }  $$
which arises from the relation  $ \, det_q \big( \rho_{ij}
\big) = 1 \, $,  for  $ \, q=1  \, $  turns into
  $$  \chi_1 + \chi_2 + \cdots + \chi_{n+1} = 0  $$
that is a relation like  $ \, Tr(x) = 0 \, $.
                                                \par
  Moreover, relations (5.2) should have quantum counterparts
  $$  "\hbox{$ q $--analogue of \ } h_i" = r_{i,i} \cdot r_{i+1,i+1}^{\, -1} \;
,  \qquad  "\hbox{$ q $--analogue of \ } \ell_i" = r_{1,1} r_{2,2} \cdots
r_{i,i} \; ,   \eqno (5.3)  $$
which should yield
  $$  {\; r_{i,i} \cdot r_{i+1,i+1}^{\, -1} - 1 \; \over \; q-1
\;}{\Big\vert}_{q=1} = h_i \; ,  \qquad  {\; r_{1,1} r_{2,2} \cdots r_{i,i} - 1
\; \over \; q-1 \;}{\Big\vert}_{q=1} = \ell_i \; ;   \eqno (5.4)  $$
now, the second relation in (5.3) is completely meaningful, so it moved us to
define  $ \, \psi \, $  by  $ \, \psi_i := {\; r_{1,1} r_{2,2} \cdots r_{i,i} -
1 \; \over \; q-1 \;} \, $;  on the other hand, the first one instead is
meaningless, for  $ r_{j,j}^{\, -1} $  does not exist; but since  $ \, r_{j,j}
\equiv 1 \; \mod\, (q-1) \, $,  we should have also
  $$  "\hbox{$ q $--analogue of \ } h_i" = "\hbox{$ q $--analogue of \ } h_i"
\cdot r_{i+1,i+1} = r_{i,i} \cdot r_{i+1,i+1}^{\, -1} \cdot r_{i+1,i+1} =
r_{i,i}  $$
 \eject
\noindent   whence the first relation in (5.4) turns into
  $$  {\; r_{i,i} - r_{i+1,i+1} \; \over \; q-1 \;}{\Big\vert}_{q=1} =
{\left( {\; r_{i,i} \cdot r_{i+1,i+1}^{\, -1} - 1 \; \over \; q-1 \;}
\cdot r_{i+1,i+1} \right)}{\Bigg\vert}_{q=1} = h_i  $$
which provides a completely meaningful expression for a (tentative)
$ q $--analogue  of  $ h_i \, $;  that's why defined  $ \, \varphi_i := {\;
r_{i,i} - r_{i+1,i+1} \; \over \; q-1 \;} \, $.  Notice that also  $ \, L_i
\, $,  resp.~$ \, K_i \, $,  is a  $ q $--analogue  of  $ \, \ell_i \, $,
resp.~$ \, h_i \, $:  this explain the relation
  $$  \xi_{\scriptscriptstyle P} \left( \psi_i \right) = \left( L_i ; \, 0 \atop
1 \right) + {\Cal O}(q-1) \; ,  \qquad \quad  \text{resp.} \quad
\xi_{\scriptscriptstyle P} \left( \varphi_i \right) = \left( K_i ; \, 0 \atop 1
\right) + {\Cal O}(q-1)  $$

\vskip5pt

   {\bf Conclusion.} \  The result of all the analysis above is that
  $$  \displaylines{
   {} \quad   \xi_{\scriptscriptstyle P} {\Big\vert}_{q=1} \! \left(
\fqPtildesln \Big/ (q-1) \, \fqPtildesln \right) =   \hfill {\ }  \cr
   {\ } \hfill   = \xi_{\scriptscriptstyle P} \! \left( \fqPtildesln \right)
\Big/ (q-1) \, \xi_{\scriptscriptstyle P} \! \left( \fqPtildesln \right) =
\gerU^{\scriptscriptstyle P} (\gerh) \Big/ (q-1) \, \gerU^{\scriptscriptstyle P}
(\gerh) \, ,  \cr
   {} \quad   \xi_{\scriptscriptstyle Q} {\Big\vert}_{q=1} \! \left(
\fqQtildesln \Big/ (q-1) \, \fqQtildesln \right) =   \hfill {\ }  \cr
   {\ } \hfill   = \xi_{\scriptscriptstyle Q} \! \left( \fqQtildesln \right)
\Big/ (q-1) \, \xi_{\scriptscriptstyle Q} \! \left( \fqQtildesln \right) =
\gerU^{\scriptscriptstyle Q} (\gerh) \Big/ (q-1) \, \gerU^{\scriptscriptstyle Q}
(\gerh) \, ,  \cr
   {} \quad   \xi_{\scriptscriptstyle P} {\Big\vert}_{q=1} \!
\left( \fqtildesln \Big/ (q-1) \, \, \fqtildesln \right) =   \hfill {\ }  \cr
   {\ } \hfill   = \xi_{\scriptscriptstyle P} \! \left( \fqtildesln \right)
\Big/ \, (q-1) \, \xi_{\scriptscriptstyle P} \! \left( \fqtildesln \right) =
\gerU^{\scriptscriptstyle P} (\gerh) \Big/ (q-1) \, \gerU^{\scriptscriptstyle P}
(\gerh) \, ,  \cr }  $$
because elements  $ F_i $,  $ \left({ L_i ; \, 0 \atop 1 }\right) $
--- resp.~$ \left({ K_i ; \, 0 \atop 1 }\right) $,  resp.~$ \left({ G_i ; \, 0
\atop 1 }\right) $  ---   and  $ E_i \, $  are enough to generate
$ \gerU^{\scriptscriptstyle P} (\gerh) \big/ (q-1) \, \gerU^{\scriptscriptstyle
P} (\gerh) $,  resp.~$ \gerU^{\scriptscriptstyle Q} (\gerh) \big/ (q-1) \,
\gerU^{\scriptscriptstyle Q} (\gerh) $,  resp.~$ \gerU^{\scriptscriptstyle P}
(\gerh) \big/ (q-1) \, \gerU^{\scriptscriptstyle P} (\gerh) $.  In particular,
in this sense we claim that  "$ \fqQtildesln $,  resp.~$ \fqPtildesln $,
is an approximation of  $ \calfQsln $,  resp.~$ \calfPsln $".
                                                      \par
   It is worth stressing that this implies that Theorem 4.1 is a direct
consequence of the specialization results about  $ \calfPsln $  and
$ \calfQsln $  proved in [Ga], \S 7 (Theorem 7.3); conversely, those results
follows from Theorem 4.1:

\vskip5pt

\proclaim{Theorem}  $ \calfQsln $  and  $ \calfPsln $  specialize to
$ \, \uh \, $  as Poisson Hopf coalgebras for  $ \, q \rightarrow 1 \, $.
The same holds for  $ \gerU^{\scriptscriptstyle Q} (\gerh) $  and
$ \gerU^{\scriptscriptstyle P} (\gerh) $  too.
\endproclaim

\demo{Proof}  To be short let us set  $ \, A{\big\vert}_{q=1} := A \big/ (q-1)
\, A \, $  for any  $ \kqqm $--algebra.  Now, we have
$ \, \fqPtildesln \subseteq \calfPsln \subseteq \gerU^{\scriptscriptstyle P}
(\gerh) \, $,  and the analysis above shows   --- through and together with that
in [Ga] ---   that  $ \, \fqPtildesln{\Big\vert}_{q=1} \cong
\gerU^{\scriptscriptstyle P} (\gerh){\Big\vert}_{q=1} \, $,  hence
  $$  \fqPtildesln{\Big\vert}_{q=1} \cong \calfPsln{\Big\vert}_{q=1} \cong
\gerU^{\scriptscriptstyle P} (\gerh){\Big\vert}_{q=1} \, ;  $$
the same holds with  $ Q $  instead of  $ P $.  Then the claim follows from
Theorem 4.1.   $ \square $
\enddemo

\vskip1,3truecm

\centerline { \bf  \S \; 6 \,  Generalization to
$ F_q[GL(n+1)] $ }

\vskip10pt

  {\bf 6.1  The quantum matrix-function algebra  $ \fqmn $.} \  We introduce the
 {\it quantum matrix-function algebra}  of order  $ n+1 $  ($ n \in \N $),  to
be called  $ \fqmn $,  as follows.
                                                      \par
  By definition,  $ \fqmn $  is the associative  $ \kq $--algebra with 1
         generated by\break
$ \{\, x_{ij} \mid i, j= 1, \dots, n+1 \,\} $  with relations
  $$  \displaylines {
   \hfill   x_{ij} x_{ik} = q \, x_{ik} x_{ij} \; ,  \qquad \quad  x_{ik} x_{hk}
= q \, x_{hk} x_{ik}   \hfill  \forall\; j<k, i<h \quad  \cr
   \hfill   x_{il} x_{jk} = x_{jk} x_{il} \; ,  \qquad \quad  x_{ik} x_{jl} -
x_{jl} x_{ik} = \left( q - \qm \right) \, x_{il} x_{jk}   \hfill  \forall\;
i<j, k<l \quad  \cr }  $$
(that is, the same of  $ \fqPsln $  but for the relation "quantum determinant =
1").  This is a bialgebra, with co-operations defined by
  $$  \Delta (x_{ij}) = \sum_{k=1}^n x_{ik} \otimes x_{kj} \, ,  \qquad
\epsilon(x_{ij}) = \delta_{ij}  \eqno  \forall\; i, j \, . \quad   $$
   \indent   From the very definition we get a bialgebra epimorphism
  $$  \pi : \fqmn \llongtwoheadrightarrow \fqPsln \, ,  \quad  x_{ij}
\mapsto \rho_{ij}   \qquad \qquad \forall \; i, j \, .  $$

\vskip7pt

  {\bf 6.2  The quantum function algebra  $ F_q[GL(n+1)] $.} \  The element
$ det_q \big( x_{ij} \big) $  of  $ \fqmn $  is group-like and central; thus by
localization at  $ det_q $  one can define a new algebra, namely
$ \fqmn \! \left[ det_q^{\, -1} \right] $:  this is now a Hopf algebra, with
bialgebra structure given by extension of that of  $ \fqmn $  and antipode
defined by
  $$  S(x_{ij}) := {(-q)}^{j-i} {det}_q \left( {(x_{hk})}_{h \neq j}^{k \neq i}
\right)   \eqno \forall\, i, j \, . \quad  $$
   \indent   The Hopf algebra  $ \, F_q[GL(n+1)] := \fqmn \left[ det_q^{\, -1}
\right] \, $  is  {\it the quantum function algebra of the group}  $ GL(n+1) $
(cf.~[Ta]).  It is clear that the bialgebra epimorphism  $ \, \pi : \fqmn
{\relbar\joinrel\twoheadrightarrow} \fqPsln \, $  extends to  $ \, \pi :
F_q[GL(n+1)] \llongtwoheadrightarrow \fqPsln \, $,  a Hopf algebra
epimorphism whose kernel is the (Hopf) ideal generated by  $ \, \Big( det_q
\big( x_{ij} \big) - 1 \Big) $.
                                                  \par
   The constructions and results in \S\S 2--4 can be easily extended to
$ F_q[GL(n+1)] $.  It is straightforward to check that  $ \kqqm $--integer
forms  $ \widetilde{F}_q^{\scriptscriptstyle P}[GL(n+1)] $  and
$ \widetilde{F}_q[GL(n+1)] $  of  $ F_q[GL(n+1)] $  can be defined mimicking
the definitions of  $ \fqPtildesln $  and  $ \fqtildesln $,  with  $ x_{ij} $'s
instead of  $ \rho_{ij} $'s:  then one has a presentation of these algebras by
generators and relations (namely, the same as for  $ \fqPtildesln $,
resp.~$ \fqtildesln $,  but for the relation  $ \, \psi_{n+1} = - \!
\sum_{\sigma \in S_{n+1} \setminus \{1\}} {(-q)}^{l(\sigma)} \,
{\left( q - 1 \right)}^{e(\sigma)-1} {\left( 1 + \qm \right)}^{e(\sigma)}
r_{1,\sigma(1)} r_{2,\sigma(2)} \cdots r_{n+1,\sigma(n+1)} $,
resp.~$ \, \sum_{i=1}^{n+1} r_{1,1} r_{2,2} \cdots r_{i-1,i-1} \chi_i =
\sum_{\sigma \in S_{n+1} \setminus \{1\}} {(-q)}^{l(\sigma)} \,
{(q-1)}^{e(\sigma)-1} {\left( 1 + \qm \right)}^{e(\sigma)}
                                    \cdot $\break
$ \cdot r_{1,\sigma(1)} r_{2,\sigma(2)} \cdots r_{n+1,\sigma(n+1)} \, $).
The upshot is that the specialization at  $ \, q=1 \, $  of both of these
integer forms is a Poisson Hopf  coalgebra isomorphic to  $ U \! \left( \gerh'
\right) $,  where  $ \gerh' $  is the Lie bialgebra obtained by central
extension of  $ \gerh $  by an element  $ \, c \, $  (namely,  $ \, c =
{\chi_{n+1}}{\big\vert}_{q=1} \, $)  such that
  $$  \displaylines{
   \qquad \qquad \qquad  c x - x c = 0  \qquad \forall\, x \in \gerh  \qquad
\quad \hbox{(i.e.  $ c $  is central)}  \cr
   \Delta(c) = c \otimes 1 + 1 \otimes c \, ,  \quad  \epsilon(c) = 0 \, ,
\quad  S(c) = -c  \cr
   \delta(c) = 4 \cdot \sum_{k=1}^n \text{f}_{n+1,k} \wedge
\text{e}_{k,n+1}  \cr }  $$
   \indent   Thus again the quantum function algebra  $ F_q[GL(n+1)] $  can be
seen as a quantum enveloping algebra, namely sort of a  "$ U_q \! \left( \gerh'
\right) $".

\vskip1,3truecm

\centerline { \bf  \S \; 7 \,  PBW theorems }

\vskip10pt

   In this section we shall prove some PBW theorems for
$ \fqPsln $:  that is, we shall exhibit some  $ \kq $--basis  of ordered
monomials in the  $ \rho_{ij} $'s  for this algebra.  To begin withs, we recall
(cf.~[Ko], [PW]) that, whenever we fix any total order in the set of generators
$ \{\, x_{ij} \mid i, j= 1, \dots, n+1 \,\} $,  the following PBW-type theorem
holds for  $ \fqmn $:

\vskip7pt

\proclaim{Proposition 7.1}  The set of ordered monomials in the generators
$ x_{ij} $'s  is
                 a  $ \kq $--basis\break
of  $ \fqmn \, $.   $ \square $
\endproclaim

\vskip7pt

  Now we wish to prove a similar result for  $ \fqPsln $;  to this end we need a
"triangularization argument", which is now explained.  Define
                                                     \par
   \indent   $ N_+ := \kq $--subalgebra  of $ \fqmn $  generated by  $ \{\,
x_{ij} \mid j<n+2-i \,\} $,
                                                     \par
   \indent   $ N_0 := \kq $--subalgebra  of $ \fqmn $  generated by  $ \{\,
x_{ij} \mid j=n+2-i \,\} $,
                                                     \par
   \indent   $ N_- := \kq $--subalgebra  of $ \fqmn $  generated by  $ \{\,
x_{ij} \mid j>n+2-i \,\} $;
                                                     \par
\noindent   then we have the following result, whose proof easily follows from
definitions and Proposition 7.1:

\vskip7pt

\proclaim{Proposition 7.2}  Let any total order in  $ \{\, x_{ij} \mid i,j = 1,
\dots, n+1 \,\} $  be fixed.  Then:
                                         \hfill\break
   \indent (a) \  the set of ordered monomials
  $$  \Big\{\, \prod_{j<n+2-i} x_{ij}^{m_{ij}} \,\Big\vert\, m_{ij} \in \N \, \;
\forall\, i,j \,\Big\}  $$
is a  $ \kq $--basis  of  $ N_+ \, $;  the set of ordered monomials
  $$  \Big\{\, \prod_i x_{i,n+2-i}^{m_i} \,\Big\vert\, m_i \in \N \, \;
\forall\, i \,\Big\}  $$
is a  $ \kq $--basis  of  $ N_0 \, $;  the set of ordered monomials
  $$  \Big\{\, \prod_{j>n+2-i} x_{ij}^{m_{ij}} \,\Big\vert\, m_{ij} \in \N \; \,
\forall\, i,j \,\Big\}  $$
is a  $ \kq $--basis  of  $ N_- \, $;
                                         \hfill\break
   \indent (b) \  $ N_0 $  is a commutative subalgebra of
$ \fqmn $;
                                         \hfill\break
 \eject
   \indent (c) \  (Triangular Decomposition) the multiplication in  $ \fqmn $
gives a  $ \kq $--vector  space isomorphism
  $$  \fqmn \cong N_+ \otimes N_0 \otimes N_- \; .  \quad \square  $$
\endproclaim

\vskip7pt

  We are now ready for the first PBW theorem for
$ \fqPsln $.

\vskip7pt

\proclaim{Theorem 7.3 ($ 1^{\text{st}} $  PBW theorem for  $ \fqPsln $)}  Let
$ \preceq $  be any fixed total ordering of the index set  $ \{\, (i,j) \mid
i,j=1,\dots,n \,\} $  such that  $ \, (i,j) \preceq (h,k) \preceq (l,m) \, $  for all  $ i $,  $ j $,  $ h $,  $ k $,  $ l $,  $ m $  such that  $ \, i < n +
1 - j \, $,  $ \, h = n + 1 - k \, $,  $ \, l > n + 1 -m \, $. Then the set of
ordered monomials
  $$  M' := \left\{\, \prod_{i<n+2-j} \!\! \rho_{ij}^{N_{ij}} \prod_{k=n+2-h}
\!\! \rho_{hk}^{N_h} \prod_{l>n+2-m} \!\! \rho_{lm}^{N_{lm}} \,\bigg\vert\,
N_{st} \in \N \; \forall \, s, t \, ; \, min\{N_1, \dots, N_{n+1}\} = 0
\,\right\}  $$
is a  $ \kq $--basis  of  $ \fqPsln \, $.
\endproclaim

\demo{Proof}  We prove now that the above set does span  $ \fqPsln $;  the
linear independence will be an easy consequence of Theorem 7.4 below.
                                                   \par
  Since  $ \fqPsln $  is a homomorphic image of  $ \fqmn $,  it is clear that
the whole set (without restriction on the indices  $ N_h $'s)  of ordered
monomials in the  $ \rho_{ij} $'s  does span  $ \fqPsln $  over  $ \kq$.  Now
pick up any monomial
  $$ \, m(\underline{N}) := \prod_{i<n+2-j}
\rho_{ij}^{N_{ij}} \cdot \prod_{k=n+2-h} \rho_{hk}^{N_h} \cdot \prod_{l>n+2-m}
\rho_{lm}^{N_{lm}} \, $$
%£
such that  $ \, d:= min\{N_1, \dots, N_{n+1}\} > 0
\, $;  since the generators  $ \rho_{h,n+2-h} $  ($ h = 1 \dots, n+1 $)  commute
with each other, we can single out of the  "$ \pi(N_0) $--part"  (with respect
to the triangular decomposition inherited from that in Proposition 7.2(c))  $ \,
n_0 := \prod_{k=n+2-h} \rho_{hk}^{N_h} = \prod_{h=1}^{n+1} \rho_{h,n+2-h}^{N_h}
\, $  of  $ m(\underline{N}) $  a factor  $ \, \rho_{1,n+1} \rho_{2,n} \cdots
\rho_{n+1,1} \, $,  and we can do it  $ d $  times.  Now using the relation
$ \, {det}_q \big( \rho_{ij} \big) - 1 = 0 \, $  we substitute the factor
$ \, \rho_{1,n+1} \rho_{2,n} \cdots \rho_{n+1,1} \, $
in  $ m(\underline{N}) $  with
  $$  {(-q)}^{{n+1 \choose 2}} - \sum_{\sigma \in S_{n+1} \setminus \{w_0\}}
{(-q)}^{{n+1 \choose 2} - \ell(\sigma)} \cdot \rho_{1,\sigma(1)}
\rho_{2,\sigma(2)} \cdots \rho_{n+1,\sigma(n+1)}   \eqno (7.1)  $$
where  $ w_0 $  is the longest element of  $ S_{n+1} \, $;  now look at the
                         various summands\break
$ \, \rho_{1,\sigma(1)} \rho_{2,\sigma(2)} \cdots \rho_{n+1,\sigma(n+1)} \, $
(up to the proper coefficient) coming in, with  $ \, \sigma \neq w_0 \, $:
whenever we have  $ \, \sigma(j+1) < \sigma(j) \, $  the commutation rules give
  $$  \rho_{j,\sigma(j)} \rho_{j+1,\sigma(j+1)} = \rho_{j+1,\sigma(j+1)}
\rho_{j,\sigma(j)} \, ;  $$
in particular, if  $ \, j \geq n + 2 - \sigma(j) \, $  and  $ \, j + 1 \leq
n +2 - \sigma(j+1) \, $  we have exactly  $ \, \rho_{j,\sigma(j)}
\rho_{j+1,\sigma(j+1)} = \rho_{j+1,\sigma(j+1)} \rho_{j,\sigma(j)} \, $;
it follows that we can factor out the monomial  $ \, \rho_{1,\sigma(1)}
\rho_{2,\sigma(2)} \cdots \rho_{n+1,\sigma(n+1)} \, $  as
  $$  \rho_{1,\sigma(1)} \rho_{2,\sigma(2)} \cdots \rho_{n+1,\sigma(n+1)} = n_+
\cdot n'_0 \cdot n_-  $$
with  $ \, n_+ \in \pi(N_+) $,  $ n'_0 \in \pi(N_0) $,  $ n_- \in \pi(N_-) $,
and  $ \, deg\left(n'_0\right) < deg(n_0) \, $  as monomials in the
$ \rho_{ij} $'s.  Therefore at the end we are left with a new expression of
$ m(\underline{N}) $  as a linear combination of monomials which, with respect
to the triangular decomposition inherited from Proposition 7.3(c), have a
"$ \pi(N_0) $--part"  of lower degree; then a simple induction argument
finishes the proof.   $ \square $
\enddemo
 \eject

\vskip7pt

  We conclude with our second PBW theorem: this is the most interesting, because
it is directly related to the classical PBW theorem attached to the natural
triangular decomposition of  $ \uh $.

\vskip7pt

\proclaim{Theorem 7.4 ($ 2^{\text{nd}} $  PBW theorem for  $ \fqPsln $)}  Let
$ \preceq $  be any fixed total ordering of the index set  $ \{\, (i,j) \mid
i,j=1,\dots,n \,\} $  such that  $ \, (i,j) \preceq (h,k) \preceq (l,m) \, $
for all  $ i $,  $ j $,  $ h $,  $ k $,  $ l $,  $ m $  such that  $ \, i > j
\, $,  $ \, h = k \, $,  $ \, l < m \, $. Then the set of ordered monomials
  $$  M := \left\{\, \prod_{i>j} \rho_{ij}^{N_{ij}} \prod_{h=k}
\rho_{hk}^{N_{hk}} \prod_{l<m} \rho_{lm}^{N_{lm}} \,\bigg\vert\, N_{st} \in \N
\; \forall \, s, t \, ; \, min \{\, N_{1,1}, \dots, N_{n+1,n+1} \,\}
= 0 \,\right\}  $$
is a  $ \kq $--basis  of  $ \fqPsln $.
\endproclaim

\demo{Proof}  As  $ \fqmn $  is clearly  $ {\N}^{{(n+1)}^2} $--graded
(by the degree in each variable), it is also  $ \N $--graded  (by the
total degree); hence  $ \fqPsln $  inherits a filtration, arising from
the filtration of  $ \fqmn $  associated to the  $ \N $--grading,  say
  $$  F_0 \subset F_1 \subset F_2 \subset \cdots \subset F_r \subset \cdots
\subset \fqPsln = \bigcup_{r=0}^{+\infty}
F_r \; ;  $$
furthermore, we have  $ \, M = \cup_{r=0}^{+\infty} M_r \, $,  with  $ \, M_r :=
M \cap F_r \, $  for all  $ \, r \in \N $.
                                             \par
  Similarly we have  $ \, M' = \cup_{r=0}^{+\infty} M'_r \, $,  with  $ \, M'_r
:= M' \cap F_r \, $  for all  $ \, r \in \N $,  where  $ M' $  is the set of
ordered monomials defined in Proposition 7.4 above: in particular,  $ M'_r $
spans  $ F_r $  over  $ \kq $.  Finally the very definitions ensure that
  $$  \# \big( M_r \big) = \# \big( M'_r \big)   \qquad \qquad  \forall \, r \in
\N \, .   \eqno (7.2)  $$
   \indent   Now consider the specialization  $ \, \fqPsln \, {\buildrel {q
\rightarrow 1} \over \llongrightarrow} \, F[SL(n+1)] \, $  and the corresponding
set  $ M^{(1)} $  of "specialized monomials", i.e.~the image of  $ M $  under
the epimorphism  $ \, \fqPsln {\relbar \joinrel \twoheadrightarrow}
F_1^{\scriptscriptstyle P} [SL(n+1)] = F[SL(n+1)] \, $:  if we prove that  $
M^{(1)} $  is a linearly independent set (over  $ k \, $), then the same will be
true for  $ M $  (over  $ \kq \, $);  in particular  $ M_r $  will be linearly
independent  ($ \, \forall r \in \N \, $),  hence it will be a  $ \kq $--basis
of  $ F_r $   ($ \, \forall r \in \N \, $),  because of (7.2), whence finally  $
M $  will be a  $ \kq $--basis  of  $ \fqPsln $.  Thus let us prove that the
set  $ M^{(1)} $  is linearly independent over  $ k $.
                                                     \par
  Assume we have in  $ \, F_1^{\scriptscriptstyle P}[SL(n+1)] = F[SL(n+1)] \, $
a relation
  $$  \sum_{m \in M^{(1)}} a_m \cdot m = 0   \eqno (7.3)  $$
for some (finitely many)  $ \, a_m \in k \setminus \{0\}
\, $;  then (7.3) lifts up to a relation in  $ k \big[
\{z_{ij}\}_{i,j=1,\dots,n+1} \big] $
  $$  \sum_{m \in M^{(1)}} a_m \cdot m(z) = b(z) \cdot \Big( det \big( z_{ij}
\big) - 1 \Big)   \eqno (7.4)  $$
for some  $ \, b(z) \in k \big[ \{z_{ij}\}_{i,j=1,\dots,n+1} \big] \, $,  the
$ m(z) $'s  having the obvious meaning.  If  $ \, b(z) = 0 \, $,  then (7.4)
gives a non-trivial algebraic relation among the  $ z_{ij} $'s,  which is
impossible; if  $ \, b(z) \neq 0 \, $,  then the right-hand-side of (7.4)
involves   --- with non-zero coefficient ---   the monomial  $ \, z_{1,1}
z_{2,2} \cdots z_{n+1,n+1} \, $  (coming out of  $ det \big( z_{ij} \big) $),
while each monomial  $ m(z) $  in the left-hand-side does not contain the
factor  $ \, z_{1,1} z_{2,2} \cdots z_{n+1,n+1} \, $;  thus again (7.4) yields a
non-trivial algebraic relation among the  $ z_{ij} $'s,  which is impossible.
The claim follows.   $ \square $
\enddemo
 \eject

\vskip1,5truecm

\Refs
\endRefs

\vskip8pt

\smallrm

[APW] \  H.~H.~Andersen, P.~Polo, Wen Kexin,  {\smallit
Representations of quantum algebras\/},  Invent.~Math.~{\smallbf 104} (1991),
1--59.

\vskip5pt

[DKP] \  C.~De Concini, V.~G.~Kac, C.~Procesi,  {\smallit
Quantum coadjoint action\/},  Jour.~Am. Math.~Soc.~{\smallbf 5} (1992), 151 --
189.

\vskip5pt

[DL] \  C.~De Concini, V.~Lyubashenko,  {\smallit Quantum
function algebra at roots of 1\/},  Adv. Math.~{\smallbf 108} (1994), 205 --
262.

\vskip5pt

[DP] \  C.~De Concini, C.~Procesi,  {\smallit Quantum groups\/},  in  L.~Boutet
de Monvel, C.~De Concini, C.~Procesi, P.~Schapira, M.~Vergne (eds.),  {\smallit
D-modules, Representation Theory, and Quantum Groups\/},  Lectures Notes in
Mathematics  {\smallbf 1565}, Springer
$ {\scriptstyle \and} $  Verlag, Berlin--Heidelberg--New York, 1993.

\vskip5pt

[Dr] \  V.~G.~Drinfeld,  {\smallit Quantum groups\/},  Proc.~ICM Berkeley 1
(1986),  789--820.

\vskip5pt

[Ga] \  F.~Gavarini,  {\smallit  Quantization of Poisson groups\/},
to appear in  Pac.~Jour.~Math.~(preprint q-alg/9604007).

\vskip5pt

[GL] \  I.~Grojnowski, G.~Lusztig,  {\smallit  On bases of irreducible
representation of quantum group  $ {GL}_n \, $\/},  in V.~Deodhar (ed.),
{\smallit  Kazhdan-Lusztig theory and related topics\/},
Cont.~Math.~{\smallbf 139} (1992), 167--174.

\vskip5pt

[Ko] \  H.~T.~Koelink,  {\smallit  On  $ \ast $--representations of the
Hopf  $ \ast $--algebra  associated with the quantum group  $ U_q(n) $\/},
Compositio Mathematicae {\smallbf 77} (1992), 199--231.

\vskip5pt

[Ji] \  M.~Jimbo,  {\smallit  A q-analogue of U(gl(N+1)), Hecke Algebras and the
Yang Baxter equation\/},  Lett.~Math. Phys.~{\smallbf 10} (1985), 63--69.

\vskip5pt

[Pa] \  P.~Papi,  {\smallit  A characterization of a good ordering in a root
system\/},  Proc.~Amer.~Math.~Soc.~{\smallbf 120} (1994), 661--665.

\vskip5pt

[PW] \  B.~Parshall, J.~Wang,  {\smallit  Quantum linear groups\/},
Mem.~Amer.~Math.~Soc.~{\smallbf 89} (1991),
$ n^{\smallit o} $  439.

\vskip5pt

[Ta] \  M.~Takeuchi,  {\smallit  Some topics on
$ {GL}_q(n) $\/},  J.~Algebra {\smallbf 147} (1992), 379--410.

\vskip0,7truecm

\enddocument